\documentclass[12pt,thmsa]{article}
\usepackage{amsmath}
\usepackage{graphicx}
\usepackage{amsfonts}
\usepackage{amssymb}

\newtheorem{Cor}{Corollary}
\newtheorem{Def}{Definition}[section]

\newtheorem{Lem}{Lemma}
\newtheorem{Pro}{Proposition}

\newtheorem{Rmq}{Remark}

\newtheorem{Theo}{Theorem}

\newcommand{\E}{\mathbb E}

\begin{document}

\begin{center}
{\Large Robust multiple-set linear canonical analysis based on minimum covariance determinant estimator}

{\Large \bigskip}

Ulrich  DEMBY BIVIGOU and Guy Martial \ NKIET

\bigskip

Universit\'{e} des Sciences et Techniques de Masuku

BP 943 \ Franceville, Gabon

E-mail : dastakreas@yahoo.fr,  gnkiet@hotmail.com.

\bigskip
\end{center}

\noindent\textbf{Abstract.} By deriving influence functions related to multiple-set linear canonical analysis (MSLCA) we show that the classical version of this analysis,  based on  empirical covariance operators, is not robust.  Then, we introduce a robust version of MSLCA by using the MCD estimator of the covariance operator of the involved random vector. The related influence functions are then derived and are shown to be bounded. Asymptotic properties of the introduced robust MSLCA are obtained and permit to propose a robust test for mutual non-correlation. This test is shown to be robust by studying the related second order influence function under the null hypothesis.

\bigskip



\section{Introduction}
\label{intro}
Many multivariate statistical methods are based  on empirical covariance operators. That is the case for multiple regression, principal components analysis, factor analysis, linear discriminant analysis, linear canonical analysis, multiple-set linear canonical analysis, and so on. However, these empirical covariance operators are known to be extremely sensitive to outliers. That is an undesirable property that makes the preceding methods themselves sensitive to outliers. For overcoming this problem, robust alternatives for these methods  have been proposed in the literature, mainly by replacing the aforementioned empirical covariance operators by robust estimators. In this vein,   robust versions of multivariate statistical methods have been introduced, especially for multiple regression (\cite{rousseuw}), principal components analysis (\cite{crouxhaesbroeckb},\cite{cui},\cite{ibazizen},\cite{serneels}), factor analysis (\cite{pison}), linear discriminant analysis (\cite{crouxfilmo},\cite{crouxjoosens},\cite{hawkins}), linear canonical analysis (\cite{crouxdehon},\cite{taskinen}). Multiple-set linear canonical analysis (MSLCA)  is an important multivariate  statistical method that analyzes  the relationship  between more than  two random vectors, so generalizing linear canonical analysis. It has been introduced for many years (e.g., \cite{gifi}) and has been studied     since then under different aspects (e.g., \cite{hwang},\cite{takane},\cite{tenenhaus}). A formulation of MSLCA within the context of Euclidean random variables has been made recently (\cite{nkiet}) and permitted to obtain an asymptotic theory for this analysis when it is estimated by using empirical covariance operators. To the best of our knowledge, such estimation of MSLCA is the one that have been tackled in the literature, despite the fact that it is known to be nonrobust as it is sensitive to outliers. So, there is a real interest in introducing a robust estimation of MSLCA as it was done for the others multivariate statistical methods. This can be done by using  robust estimators  of the covariance operators  of the involved random vectors instead of the empirical covariance operators. Among such robust estimators, the minimum covariance determinant (MCD) estimator has been extensively studied (\cite{butler}, \cite{cator},\cite{cator2},\cite{crouxhaesbroecka}), and it is known to have good  robustness properties. Also,  its asymptotic properties have been obtained (\cite{butler},\cite{cator},\cite{cator2}) mainly under elliptical distribution.

In this paper, we propose a robust version of MSLCA based on MCD estimator of the covariance operator. We start by recalling, in Section 2, the notion of MSLCA for Euclidean random variables and we study its robustness properties by deriving  the  influence functions of  the functionals that lead  to its  estimator  from  the empirical covariance operators. It is  proved that the influence function of the operator that determines MSLCA is not bounded. In Section 3, we introduce a robust estimation of MSLCA (denoted by RMSLCA)  by using the MCD estimator of the covariance operator on which this analysis is defined. Then we derive the influence function of the  operator that determines RMSLCA, which is proved to be bounded, and that  of the canonical coefficients and the canonical directions. Section 4  is devoted to asymptotic properties of  RMSLCA. We obtain limiting distributions that are then used in Section 5 where a robust test for mutual non-correlation is introduced. The robustness properties of this test are studied through the derivation of the second order influence function of the test statistic under the null hypothesis. The proofs of all theorems and propositions are postponed in Section 6.

\section{Influence in multiple-set canonical analysis}
In this section we  recall the notion of multiple-set linear canonical analysis (MSLCA) of Euclidean random variables as introduced by Nkiet\cite{nkiet},  and also its   estimation based on empirical covariance operators. Then, the robustness properties of this analysis are  studied through derivation of the  influence functions that correspond to the   functionals  related to it.
\subsection{Multiple-set linear canonical analysis}
\noindent  Letting  $(\Omega,\mathcal{A},P) $ be a probability space, and  $K$ be an integer such that $K\geq 2$, we consider  random variables  $X_1,\cdots,X_K$ defined on this  probability space and with values  in   Euclidean vector spaces $\mathcal{X}_1,\cdots,\mathcal{X}_K$ respectively. We then consider the space  $\mathcal{X}:=\mathcal{X}_1\times\cdots\mathcal{X}_K$ which is also an Euclidean vector space equipped with the inner product $\langle\cdot,\cdot\rangle_\mathcal{X}$ defined by:
\[
\forall\alpha\in \mathcal{X},\,\forall\beta\in\mathcal{X},\,\,\langle\alpha,\beta\rangle_\mathcal{X}=\sum_{k=1}^K\langle\alpha_k,\beta_k\rangle_k,
\]
where   $\langle\cdot ,\cdot\rangle _ k$ is the inner product of $\mathcal{X}_k$ and $\alpha=\left(\alpha_1,\cdots,\alpha_K\right)$, 
$\beta=\left(\beta_1,\cdots,\beta_K\right)$. 
From now on, we assume that the following assumption holds :

\bigskip

\noindent ($\mathcal{A}_1$):  for   $k\in\{1,\cdots, K\}$, we have $\mathbb{E}(X_k)=0$  and $\mathbb{E}(\Vert X_k\Vert _k^2)<+\infty $, where $\Vert\cdot\Vert _k$ denotes the norm induced by   $\langle\cdot ,\cdot\rangle _ k$ .  

\bigskip

\noindent Then, we consider the  random vector
$
X=\left(X_1,\cdots,X_K\right)
$
with values in  $\mathcal{X}$, and we can give the following definition    of multiple-set  linear canonical analysis  (see  \cite{nkiet}):

\bigskip

\begin{Def}\label{mslca}
The multiple-set linear canonical analysis
(MSLCA) of  $X$ is the search of a sequence  $\left(\alpha^{(j)}\right)_{1\leq j\leq q}$ of vectors of  $\mathcal{X}$, where $q=\dim(\mathcal{X})$,  satisfying:
\begin{equation}\label{sol1}
\alpha^{(j)}=\arg\max_{\alpha\in C_j}\mathbb{E}\left(<\alpha , X>_ \mathcal{X}^2\right),
\end{equation}
where
$
C_1=\left\{\alpha\in \mathcal{X}\,/\,\sum_{k=1}^Kvar(<\alpha_k,X_k>_k)=1\right\}$
and for  $j\geq 2$ :
\begin{equation*}\label{contr2}
C_j=\left\{\alpha\in  C_1\,/\,\sum_{k=1}^Kcov\left(<\alpha^{(r)}_k,X_k>_k,<\alpha_k,X_k>_k\right)=0,\,\,\forall r\in\{1,\cdots,j-1\}\right\}.
\end{equation*}
\end{Def}

\bigskip

\noindent A solution of the above maximization  problem is obtained from spectral analysis of an  operator  that will know be specified. For $(k,\ell)\in\{1,\cdots,K\}^2$, let us consider the covariance operators 
\[
V_{k\ell}=\mathbb{E}\left(X_\ell\otimes X_k\right)=V_{\ell k}^\ast \textrm{ and } V_k:=V_{kk},
\]
where   $\otimes $ denotes the tensor product such that  $x\otimes y $ is the linear map : $h\mapsto <x,h>y$, and $A^\ast $ denotes the adjoint of $A$. Letting  $\tau_k$ be the canonical projection 
\[
\tau_k\,:\,\alpha\in \mathcal{X}\mapsto \alpha_k\in \mathcal{X}_k
\]
which  adjoint operator  $\tau_k^\ast$   is the map
\[
\tau_k^\ast\,:\,t\in \mathcal{X}_k\mapsto (\underbrace{0,\cdots,0}_{k-1\textrm{ times}},t,0,\cdots , 0)\in \mathcal{X},
\]
 we  consider the operators  defined as
\begin{equation}\label{phipsi}
\Phi=\sum_{k=1}^K\tau_k^\ast V_k\tau_k\,\,\,\,
\textrm{ and }\,\,\,\,
\Psi=\sum_{k=1}^K\sum_{\stackrel{\ell=1}{\ell\neq k}}^K\tau_k^\ast V_{k\ell}\tau_\ell.
\end{equation}
The covariance operator $V_k$ is a self-adjoint and positive operator; we assume throughout this paper  that it is invertible. Then, it is easy to check that $\Phi $ is also self-adjoint positive and invertible  operator, and we consider 
\[
T=\Phi^{-1/2}\Psi\Phi^{-1/2}. 
\]
The spectral analysis of this last operator gives a solution of the maximization problem specified in Definition \ref{mslca}. Indeed, if $\left\{
\beta^{(1)},\cdots,\beta^{(q)}\right\} $  is an orthonormal basis of   $\mathcal{X}$  such that  $\beta^{(j)}$ is an eigenvector of   $T$ associated with the  $j$-th largest eigenvalue  $\rho_j$, then we obtain a solution  of (\ref{sol1}) by taking 
$\alpha^{(j)}=\Phi^{-1/2}\beta^{(j)}$, 
 and we have  $\rho_j=<\beta^{(j)},T\beta^{(j)}>_\mathcal{X}$. Finally, the MSLCA of $X$ is  the family $\left(\rho_j,\alpha^{(j)}\right)_{1\leq j\leq q}$ obtained as indicated above. The $\rho_j$'s are termed  \textit{the canonical coefficients} and the $\alpha^{(j)}$'s are termed \textit{the  canonical directions}. 

\bigskip

\noindent Note that $T$ can be expressed as a function of the covariance operator $V=\mathbb{E}\left(X\otimes X\right)$  of $X$. Indeed, denoting by  $\mathcal{L}(\mathcal{X})$ the space of linear maps fom $\mathcal{X}$ to itself, and considering the linear maps $f$ and $g$ from   $\mathcal{L}(\mathcal{X})$   to itself defined as
\begin{equation}\label{fg}
f(A)=\sum_{k=1}^K\tau_k^\ast \tau_kA\tau_k^\ast\tau_k\,\,\,\,
\textrm{ and }\,\,\,\,
g(A)=\sum_{k=1}^K\sum_{\stackrel{\ell=1}{\ell\neq k}}^K\tau_k^\ast\tau_kA\tau_\ell^\ast\tau_\ell,
\end{equation}
it is easy to check, by using properties of tensor produts (see \cite{dauxois}), that
\begin{equation}\label{opcov}
V_{k\ell}=\mathbb{E}\bigg((\tau_\ell (X))  \otimes (\tau_k (X))\bigg)=\tau_k\mathbb{E}\left(X\otimes X\right)\tau_\ell^\ast=\tau_kV\tau_\ell^\ast,\,\,\,\,\,V_{k}=\tau_kV\tau_k^\ast
\end{equation}
and, therefore, from (\ref{phipsi}), (\ref{fg}) and (\ref{opcov}), it follows
\begin{equation*}
T= f\left(V\right)^{-1/2}g\left(V\right)f\left(V\right)^{-1/2}.
\end{equation*}

\subsection{Estimation based on empirical covariance operator}
\noindent Now, we recall the  classical way for estimating  MSLCA by using empirical covariance operators (see, e.g., \cite{nkiet}). 
For $k=1,\cdots, K$, let $\{X_k^{(i)}\}_{1\leq i\leq n}$  be  an i.i.d. sample  of $X_k$. We then  consider the sample means and empirical covariance operators defined for $(k,\ell)\in\{1,\cdots,K\}^2$ as
\begin{eqnarray*}
\overline{X}_{k\cdot n}=\frac{1}{n}\sum_{i=1}^nX^{(i)}_k,\,\,\,\,\,\,
\widehat{V}_{k\ell\cdot n}&=&\frac{1}{n}\sum_{i=1}^n\left(X^{(i)}_\ell-\overline{X}_{\ell\cdot n}\right)\otimes \left(X^{(i)}_k-\overline{X}_{k\cdot n}\right)\\
&=& \frac{1}{n}\sum_{i=1}^nX^{(i)}_\ell\otimes X^{(i)}_k-\overline{X}_{\ell\cdot n}\otimes\overline{X}_{k\cdot n},
\end{eqnarray*}
and $\widehat{V}_{k\cdot n}:=\widehat{V}_{kk\cdot n}$. These permit to define  random operators, with values in    $\mathcal{L}(\mathcal{X})$,  as
\begin{equation}\label{phipsin}
\widehat{\Phi}_n=\sum_{k=1}^K\tau_k^\ast \widehat{V}_{k\cdot n}\tau_k\,\,\,\,
\textrm{ and }\,\,\,\,
\widehat{\Psi}_n=\sum_{k=1}^K\sum_{\stackrel{\ell=1}{\ell\neq k}}^K\tau_k^\ast \widehat{V}_{k\ell\cdot n}\tau_\ell
\end{equation}
and to  estimate $T$ by  
\begin{equation}\label{estimt}
\widehat{T}_n=\widehat{\Phi}_n^{-1/2}\widehat{\Psi}_n\widehat{\Phi}_n^{-1/2}.
\end{equation}
Considering the eigenvalues $\widehat{\rho}_{1\cdot n}\geq\widehat{\rho}_{2\cdot n}\cdots\geq\widehat{\rho}_{q\cdot n}$ of  $\widehat{T}_n$, 
and  $\left\{
\widehat{\beta}^{(1)}_n,\cdots,\widehat{\beta}^{(q)}_n\right\} $   
an orthonormal basis of   $\mathcal{X}$  such that  $\widehat{\beta}^{(j)}_n$ is an eigenvector of $\widehat{T}_n $ associated with $\widehat{\rho}_{j\cdot n}$, we can estimate $\rho_j$ by $\widehat{\rho}_{j\cdot n}$,  $\beta^{(j)}$ by $ \widehat{\beta}^{(j)}_n$ and $\alpha^{(j)}$ by $ \widehat{\alpha}^{(j)}_n=\widehat{\Phi}_n^{-1/2}\widehat{\beta}^{(j)}_n$. 

\bigskip

\noindent The random operator $\widehat{T}_n$ can also be expressed as a function of the empirical covariance operator of the $X^{(i)}$'s  that are defined as 
$
X^{(i)}=\left(X^{(i)}_1,\cdots,X^{(i)}_K\right)$; this
empirical covariance operator is 
\begin{eqnarray*}
\widehat{V}_{ n}&=&\frac{1}{n}\sum_{i=1}^n\left(X^{(i)}-\overline{X}_{ n}\right)\otimes \left(X^{(i)}-\overline{X}_{n}\right)\nonumber\\
&=& \frac{1}{n}\sum_{i=1}^nX^{(i)}\otimes X^{(i)}-\overline{X}_{n}\otimes\overline{X}_{ n}.
\end{eqnarray*}
where
$
\overline{X}_{n}=n^{-1}\sum_{i=1}^nX^{(i)}.
$
Since   $\widehat{V}_{k\ell \cdot n}=\tau_k\widehat{V}_{ n}\tau_\ell^\ast$ and  $\widehat{V}_{k\cdot n}=\tau_k\widehat{V}_{ n}\tau_k^\ast$, we  straighforwardly obtain from   (\ref{fg}), (\ref{phipsin}) and (\ref{estimt}):
\begin{equation}\label{tnfg}
\widehat{T}_n= f\left(\widehat{V}_{ n}\right)^{-1/2}g\left(\widehat{V}_{ n}\right)f\left(\widehat{V}_{ n}\right)^{-1/2}.
\end{equation}

\subsection{Influence functions}
For studying the effect of a small amount of  contamination at a given point    on  MSLCA it is important,  as usual in robustness litterature (see \cite{hampel}), to use influence function. More precisely, we have to derive expressions of the  influence functions related to the functionals that give $T$, $\rho_j$ and $\alpha^{(j)}$ (for $1\leq j \leq q$) at  the distribution $\mathbb{P}_X$ of $X$. Recall that the influence function of a functional $S$ at $\mathbb{P}$ is defined as
\[
\textrm{IF}\left(x;S,\mathbb{P}\right)=\lim_{\varepsilon\downarrow 0}\frac{S\left((1-\varepsilon)\mathbb{P}+\varepsilon\delta_x\right)-S(\mathbb{P})}{\varepsilon},
\]
where $\delta_x$ is the Dirac measure putting all its mass  in $x$.

\bigskip

\noindent First, we  have to specify the functionals related to $T$,   $\rho_j$ and $\alpha^{(j)}$ (for $1\leq j \leq q$) and their empirical counterparts. Let us consider the functional $\mathbb{T}$ given by 
\[
\mathbb{T}(\mathbb{P})= f\left(\mathbb{V}(\mathbb{P} )\right)^{-1/2}g\left(\mathbb{V}(\mathbb{P} )\right)f\left(\mathbb{V}(\mathbb{P} )\right)^{-1/2}.
\]
where $\mathbb{V}$ is the functional defined as
\[
\mathbb{V}(\mathbb{P} )=\int x\otimes x\,d\mathbb{P}(x)-\left(\int x\,d\mathbb{P}(x)\right)\otimes \left(\int x\,d\mathbb{P}(x)\right).
\]
Applying this functional to the distribution $\mathbb{P}_X$ of $X$ gives  $\mathbb{V}(\mathbb{P}_X)=V$ and, therefore,   $\mathbb{T}(\mathbb{P}_X)=T$. For $j\in\{1,\cdots,q\}$, denoting by $\mathbb{R}_j$ (resp. $\mathbb{B}_j$; resp. $\mathbb{A}_j$ ) the functional such that $\mathbb{R}_j(\mathbb{P})$ is the $j$-th largest eigenvalue of $\mathbb{T}(\mathbb{P})$ (resp. the associated eigenvector; resp. $\mathbb{A}_j(\mathbb{P})= f\left(\mathbb{V}(\mathbb{P} )\right)^{-1/2}\mathbb{B}_j(\mathbb{P})$ ), we have $\mathbb{R}_j(\mathbb{P}_X)=\rho_j$,  $\mathbb{B}_j(\mathbb{P}_X)=\beta^{(j)}$  and $\mathbb{A}_j(\mathbb{P}_X)=\alpha^{(j)}$. 

\bigskip

\noindent Furthermore, denoting by $\mathbb{P}_{ n}$ the empirical measure corresponding to the sample $\{X^{(1)},\cdots,X^{(n)}\}$, we have  
\[
\mathbb{V}(\mathbb{P}_{ n})=\widehat{V}_n,   \,\,\,\mathbb{T}(\mathbb{P}_{n})=\widehat{T}_n, \,\,\,\mathbb{R}_j(\mathbb{P}_{n})=\widehat{\rho}_{j\cdot n}, \,\,\,\mathbb{B}_j(\mathbb{P}_{n})=\widehat{\beta}^{(j)}_n \,\,\,\textrm{ and } \,\,\,\mathbb{A}_j(\mathbb{P}_{ n})=\widehat{\alpha}^{(j)}_n.
\]
These functionals are to be taken into account in order to derive the influence functions related to MSLCA. We make the following assumption:

\bigskip

\noindent  $(\mathcal{A}_2)$  : For all $ k \in \left\{1, ..., K\right\}$, we have    $V_k = I_k$, where \ $I_k$\  denotes the identity operator of  $ \mathcal{X}_k$. 

\bigskip

\noindent Then, we have the following theorem that gives the influence function of $T$.

\begin{Theo}\label{ift}
We suppose that  the assumptions  $(\mathcal{A}_1)$ and $(\mathcal{A}_2)$ hold. Then,    for any vector $x=(x_1,\cdots,x_K)\in\mathcal{X}$ we have:
\begin{equation}\label{if1}
\textrm{IF}(x;T ,\mathbb{P}_X ) =\sum_{k=1}^K \sum_{\underset{\ell\neq k}{\ell=1}}^K  -\frac{1}{2}\tau_k^*\left(x_k\otimes x_k\right) V_{k\ell}\tau_\ell - \frac{1}{2} \tau_\ell^* V_{\ell k}\left(x_k\otimes x_k\right)\tau_k+\tau_k^* \left(x_\ell\otimes x_k\right)\tau_\ell.
\end{equation}
\end{Theo}

\noindent As $T$ determines MSLCA, it is important to ask whether its influence  function is bounded. If so, we say that  MSLCA  is robust  because it would mean that a contamination at the point $x$ has a limited effect on $T$.  The following proposition shows that $\textrm{IF}(x;T,\mathbb{P}_X)$ is not bounded. We denote by $\Vert\cdot\Vert_{\mathcal{L}(\mathcal{X})}$ the operators norm defined as  $\Vert A\Vert_{\mathcal{L}(\mathcal{X})}=\sqrt{\textrm{tr}\left(AA^\ast\right)}$.

\bigskip

\begin{Pro}\label{nonbornee}
We suppose that  the assumptions  $(\mathcal{A}_1)$ and   $(\mathcal{A}_2)$ hold. Then,  there exists $x_0\in\mathcal{X}$ such that: 
\[
\lim_{t\rightarrow +\infty}\Vert\textrm{IF}(tx_0;T,\mathbb{P}_X)\Vert_{\mathcal{L}(\mathcal{X})}=+\infty.
\]
\end{Pro}
\noindent Now, we give in the following theorem, the influence functions related to the canonical coefficients and the   canonical directions.

\begin{Theo}\label{ifcano1}
We suppose that  the assumptions  $(\mathcal{A}_1)$ and   $(\mathcal{A}_2)$ hold. Then,  for any $x\in\mathcal{X}$ and any $j\in\{1,\cdots,q\}$, we have:
\begin{eqnarray*}
\textrm{(i) }\,\,\textrm{IF}(x;\rho_j ,\mathbb{P}_X )& =&
\sum_{k=1}^{K} \sum_{\underset{\ell\neq k}{\ell=1}}^{K} <\beta_k^{(j)} , x_k >_k< x_\ell -V_{\ell k}x_k, \beta_\ell^{(j)}>_\ell .
\end{eqnarray*}

\noindent\textrm{(ii) }\,\, We suppose, in addition, that    $\rho_1>\rho_2>\cdots >\rho_q$.  Then :
\begin{eqnarray}\label{ifalphaj}
\textrm{IF}( x , \alpha^{(j)} ,\mathbb{P}_X) 
&=&\sum_{k=1}^{K} \sum_{\underset{l\neq k}{\ell=1}}^{K}\sum_{\underset{m\neq j}{m=1}}^{q} \frac{1}{\rho_j - \rho_m}\bigg(<\beta_k^{(m)} ,  x_k  >_k< x_\ell , \beta_\ell^{(j)}  >_\ell \nonumber\\ 
& &-\frac{1}{2}  < \beta_k^{(m)} ,  x_k >_k < x_k , V_{k\ell} \beta_\ell^{(j)}  >_k \nonumber\\  
& &- \frac{1}{2} < x_k ,  V_{k\ell }\beta_\ell^{(m)}  >_\ell < x_k , \beta_k^{(j)} >_k \bigg)\beta^{(m)}\\
& &-\frac{1}{2}\left( \sum_{k=1}^{K}\left[ \tau_k^*\left(x_k\otimes x_k\right)\tau_k + < \beta_k^{(j)}, x_k>_k^2\mathbb{I}\right]  -  2\mathbb{I}\right)\beta^{(j)} \nonumber,
\end{eqnarray}
where $\mathbb{I}$ denotes the identity operator of $\mathcal{X}$.
\end{Theo}

\bigskip

\begin{Rmq}\label{rmq1}

Romanazzi\cite{romanazzi} derived  influence functions for the squared canonical coefficients and the canonical directions obtained from  linear canonical analysis (LCA) of two random vectors. LCA is in fact a particular  case of MSLCA obtained  when  $K=2$ (see \cite{nkiet}). With Theorem \ref{ifcano1} we recover the   results of \cite{romanazzi} when whe take   $K=2$. We will only show it below for the canonical coefficients.  For $j\in \left\{1, ..., q\right\}$,  by  applying Theorem \ref{ifcano1} with   $K = 2 $,  we obtain
\begin{eqnarray}\label{ifacl1}
\textrm{IF}(X;\rho_{j}, \mathbb{P}_X )&=& < \beta_1^{(j)} ,  x_1  >_1 < \beta_2^{(j)} , x_2 - V_{21}x_1>_2 \nonumber\\
& &+ < \beta_2^{(j)} ,  x_2 >_2 < \beta_1^{(j)} , x_1 - V_{12}x_2>_1\nonumber\\ 
&=&2< \beta_1^{(j)} , x_1 >_1\, < x_2 , \beta_2^{(j)} >_2  - < \beta_1^{(j)} , x_1  >_1\, <x_1 , V_{12} \beta_2^{(j)}>_1\nonumber\\
&  & -<x_2 , V_{21}\beta_1^{(j)}>_2 < x_2 , \beta_2^{(j)}>_2.
\end{eqnarray}
The linear canonical analysis (LCA) of $X_1$ and $X_2$ is obtained from the spectral analysis of $R=V_{12}V_{21}$  (since $V_1=I_1$ and $V_2=I_2$). If we denote by $\lambda_j$, $\eta_1^{(j)}$, $\eta_2^{(j)}$  the related squared canonical coefficients and canonical vectors, it is known (see Remark 2.2 in  \cite{nkiet}) that 
\begin{equation}\label{relcano}
\lambda_j=\rho_j^2,\,\,\, \eta_k^{(j)}=\sqrt{2} \beta_k^{(j)}\,\,\,  (k=1,2),\,\,\, V_{12}\eta_2^{(j)}=\rho_j\eta_1^{(j)},\,\,\,\textrm{  and }\,\,\, V_{21}\eta_1^{(j)}=\rho_j\eta_2^{(j)}.
\end{equation}
Then, putting  $u_j=<x_1,\eta_1^{(j)}>_1$  and  $v_j=<x_2,\eta_2^{(j)}>_2$, we deduce from (\ref{ifacl1}), (\ref{relcano})   and  the equality  $\textrm{IF}(X;\rho_{j}^2, \mathbb{P}_X )=2\rho_j\,\textrm{IF}(X;\rho_{j}, \mathbb{P}_X )$ that:
\begin{equation}\label{ifrho2}
\textrm{IF}(X;\rho_{j}^2, \mathbb{P}_X )=2\rho_j\,u_j\,v_j-\rho_j^2u_j^2-\rho_j^2v_j^2,
\end{equation}
what is the result obtained in \cite{romanazzi}.
\end{Rmq}

\section{Robust multiple-set linear canonical analysis  (RMSLCA)}

\noindent It has been seen that the MSLCA based on empirical covariance operator is not robust   since $\textrm{IF}(x;T ,\mathbb{P}_X ) $  is not bounded. There is therefore an interest in proposing a robust version of MSLCA. In this section, we introduce  such a version by replacing in (\ref{tnfg}) the empirical  covariance operator by  a  robust estimator of $V$. More precisely, we use the minimum covariance determinant (MCD) estimator of $V$. We consider the following assumption:

\bigskip

\noindent ($\mathcal{A}_3$) : the distribution $\mathbb{P}_X$ of $X$  is an  elliptical contoured distribution with density 
\[
f_X(x)=(\textrm{det}(V))^{-1/2}h(<x,V^{-1}x>_\mathcal{X}),
\]
where $h\,:\,[0,+\infty[\rightarrow  [0,+\infty[$ is  a  function having a strictly negative derivative $h^\prime$.
\bigskip

\noindent We first define the estimator of MSLCA based on MCD estimator of $V$, then we derive the  related influence functions.
\subsection{Estimation of  MSLCA based on MCD estimator}
Letting   $\gamma$ be a fixed real  such that $0 < \gamma < 1$,  we  consider a subsample $\mathcal{S} \subset \left\{X^{(1)}, ... , X^{(n)}\right\}$ of size $ h_n\geq \lceil n\gamma\rceil$, where $X^{(i)}=\left(X^{(i)}_1,\cdots,X^{(i)}_K\right)$,  and we define the empirical mean  and  covariance operator based on this subsample by: 
\[
\widehat{M}_n(\mathcal{S}) = \frac{1}{h_n}\sum_{X^{(i)} \in\mathcal{S}}X^{(i)} 
\]
and
\[
\widehat{V}_n(\mathcal{S}) = \frac{1}{h_n}\sum_{X^{(i)} \in \mathcal{S}}\left(X^{(i)} - \widehat{M}_n(\mathcal{S})\right)\otimes\left( X^{(i)} - \widehat{M}_n(\mathcal{S})\right).
\]
We denote by $\widehat{\mathcal{S}}_n$ the subsample of $\left\{X^{(1)}, ... , X^{(n)}\right\}$ which minimizes the determinant $\textrm{det}\left( \widehat{V}_n(\mathcal{S})\right)$ of $\widehat{V}_n(\mathcal{S}) $ over all subsamples of size $h_n$.
Then, the  MCD estimators of the mean  and   the covariance operator    of $X$ are $\widehat{M}_n(\widehat{\mathcal{S}}_n)$ and  $\widehat{V}_n(\widehat{\mathcal{S}}_n)$, respectively.  It is well known that the these  estimators are robusts and have high breakdown points  (see, e.g.,  \cite{rousseuw}). From them, we can  introduce an estimator of MSLCA which is expected to be also robust. Indeed, putting 
\[
\widetilde{V}_n:=\widehat{V}_n(\widehat{\mathcal{S}}_n),
\] 
we consider the random operators with values in $\mathcal{L}(\mathcal{X})$ defined as
\[
\widetilde{\Phi}_n=\sum_{k=1}^K\tau_k^\ast \widetilde{V}_{k\cdot n}\tau_k\,\,\,\,
\textrm{ and }\,\,\,\,
\widetilde{\Psi}_n=\sum_{k=1}^K\sum_{\stackrel{\ell=1}{\ell\neq k}}^K\tau_k^\ast \widetilde{V}_{k\ell\cdot n}\tau_\ell,
\]
where $\widetilde{V}_{k\cdot n}=\tau_k\widetilde{V}_{n}\tau_k^\ast$ and   $\widetilde{V}_{k\ell\cdot n}=\tau_k\widetilde{V}_{n}\tau_\ell^\ast$, and we estimate $T$ by  
\begin{equation}\label{tntilde}
\widetilde{T}_n=\widetilde{\Phi}_n^{-1/2}\widetilde{\Psi}_n\widetilde{\Phi}_n^{-1/2}.
\end{equation}
Considering the eigenvalues $\widetilde{\rho}_{1\cdot n}\geq\widetilde{\rho}_{2\cdot n}\cdots\geq\widetilde{\rho}_{q\cdot n}$ of  $\widetilde{T}_n$, 
and  $\left\{
\widetilde{\beta}^{(1)}_n,\cdots,\widetilde{\beta}^{(q)}_n\right\} $   
an orthonormal basis of   $\mathcal{X}$  such that  $\widetilde{\beta}^{(j)}_n$ is an eigenvector of $\widetilde{T}_n $ associated with $\widetilde{\rho}_{j\cdot n}$, we estimate   $\rho_j$ by $\widetilde{\rho}_{j\cdot n}$,  $\beta^{(j)}$ by $ \widetilde{\beta}^{(j)}_n$  and $\alpha^{(j)}$ by $ \widetilde{\alpha}^{(j)}_n=\widetilde{\Phi}_n^{-1/2}\widetilde{\beta}^{(j)}_n$.
This gives a robust MSLCA that we denote by RMSLCA.
\subsection{Influence functions}
In order to derive the influence functions related to the above estimator of MSLCA, we  have to specify the functional that corresponds to it. For doing that, we will first recall the functional associated to the above MCD estimator of covariance operator. Let
\[
E_\gamma=\left\{x \in  \mathcal{X} , <x ,V^{-1}x >_\mathcal{X} \,\leq r^2(\gamma)\right\}
\]
where $r(\gamma)$ is determined by the equation
\[
\frac{2\pi^{q/2}}{\Gamma (q/2)}  \int_0^{r(\gamma)}\,t^{q-1}\,h(t^2)\,dt=\gamma,
\]
$\Gamma$ being the usual gamma function. The  functional    $\mathbb{V}_{1,\gamma}$ related to the aforementioned MCD estimator of $V$    is defined in \cite{cator} (see also \cite{butler}, \cite{crouxhaesbroecka})  by
\[
\mathbb{V}_{1,\gamma}(\mathbb{P})=\frac{1}{\gamma}\int_{E_{\gamma}}\left(x - M_\mathbb{P}(E_{\gamma})\right)\otimes \left(x - M_\mathbb{P}(E_{\gamma})\right)d \mathbb{P}\left(x\right),
\]
where
\[
M_\mathbb{P}(B)=\frac{1}{\gamma}\int_{B}x \,d \mathbb{P}\left(x\right).
\]
It is known that      $\mathbb{V}_{1,\gamma}(\mathbb{P}_X)=\sigma_\gamma^2\,V$  where
\[
\sigma^2_\gamma=\frac{2\pi^{q/2}}{\gamma\,q\,\Gamma(q/2)}\int_0^{r(\gamma)}t^{q+1}h(t^2)\,dt.
\]
 Therefore,  the functional $\mathbb{T}_{1,\gamma}$ related to $T$ is defined as
\[
\mathbb{T}_{1,\gamma}(\mathbb{P})=f(\mathbb{V}_{1,\gamma}(\mathbb{P}))^{-1/2}g(\mathbb{V}_{1,\gamma}(\mathbb{P}))f(\mathbb{V}_{1,\gamma}(\mathbb{P}))^{-1/2}
\]
where $f$ and $g$ are defined in (\ref{fg}). Now, we can give the influence functions related to RMSLCA of $X$. First, putting 
\[
\kappa_0=\frac{\pi^{q/2}}{(q+2)\Gamma(q/2+1)}\int_0^{r(\gamma)}t^{q+3}\,h^\prime(t^2)\,\,dt,
\]
and $T_\gamma=\mathbb{T}_{1,\gamma}(\mathbb{P}_X)$, we have:

\bigskip

\begin{Theo}\label{iftmcd}
We suppose that the assumptions  $(\mathcal{A}_1)$ to  $(\mathcal{A}_3)$ hold. Then
\begin{eqnarray*}
\textrm{IF}\left(x;T_\gamma,\mathbb{P}_X \right)
&=&-\frac{\sigma_\gamma^{-2}}{2\kappa_0}\textrm{\large \textbf{1}}_{E_\gamma}(x)\,\,\textrm{IF}( x ; T , \mathbb{P}_X),
\end{eqnarray*}
where $\textrm{IF}( x ; T , \mathbb{P}_X)$ is given in (\ref{if1}).
\end{Theo}

\bigskip

\noindent From this theorem we  to obtain the following proposition which proves that RMSLCA is robust since the preceding influence function is bounded. We denote by $\Vert\cdot\Vert_\infty$ the usual operators  norm defined by $\Vert A\Vert_\infty=\sup_{x\in\mathcal{X}-\{0\}}\left(\Vert Ax\Vert_\mathcal{X}/\Vert x\Vert_\mathcal{X}\right)$.

\bigskip

\begin{Pro}\label{iftmcdbornee}
We suppose that the assumptions  $(\mathcal{A}_1)$ to  $(\mathcal{A}_3)$ hold.  Then,
\[
\sup_{x\in\mathcal{X}}\Vert \textrm{IF}(x;T _{\gamma}, \mathbb{P}_X)\Vert_{\infty} \leq\frac{\sigma_\gamma^{-2}}{2\vert\kappa_0\vert}K(K-1)\bigg(\Vert V\Vert_\infty+1\bigg)\Vert V^{1/2}\Vert_\infty^2\,r^2(\gamma).
\]
\end{Pro}

\noindent
\noindent Now, we give in the following theorem, the influence functions related to the canonical coefficients and the canonical directions obtained from RMSLCA. For $j\in\{1,\cdots,q\}$, denoting by $\mathbb{R}_{\gamma\cdot j}$ (resp. $\mathbb{B}_{\gamma\cdot j}$; resp. $\mathbb{A}_{\gamma\cdot j}$ ) the functional such that $\mathbb{R}_{\gamma\cdot j}(\mathbb{P})$ is the $j$-th largest eigenvalue of $\mathbb{T}_{1,\gamma}(\mathbb{P})$ (resp. the associated eigenvector; resp. $\mathbb{A}_{\gamma\cdot j}(\mathbb{P})= f\left(\mathbb{V}_{1,\gamma}(\mathbb{P} )\right)^{-1/2}\mathbb{B}_{\gamma\cdot j}(\mathbb{P})$ ), we put  $\rho_{\gamma\cdot j}=\mathbb{R}_{\gamma\cdot j}(\mathbb{P}_X)$,  $\beta^{(j)}_{\gamma}=\mathbb{B}_j(\mathbb{P}_X)$  and $\alpha^{(j)}_{\gamma}=\mathbb{A}_{\gamma\cdot j}(\mathbb{P}_X)$.  Considering
\begin{equation*}\label{const1}
\nu_0=\frac{2\pi^{q/2}}{\Gamma(q/2)}h\left(r(\gamma)^2\right)r(\gamma)^{q-1}\sigma_\gamma,\,\,\,\nu_1=\frac{r(\gamma)}{\sigma_\gamma},\,\,\,\nu_2=\frac{2\nu_0\nu_1^3}{\sigma_\gamma\, q(q+2)}-\frac{2\gamma}{\sigma_\gamma},
\end{equation*}
\begin{equation}\label{const2}
\kappa_1=-\frac{r(\gamma)^2}{q\gamma},\,\,\,\kappa_2=\frac{\sigma_\gamma\,\nu_2+2\gamma}{q\,\gamma\, \sigma_\gamma\nu_2},\,\,\,\kappa_3=-\frac{2}{\sigma_\gamma\nu_2},\,\,\,\kappa_4=\frac{r(\gamma)^2-q\sigma_\gamma^2}{q},
\end{equation}
we have:

\bigskip

\begin{Theo}\label{ifcanomcd}
We suppose that  the assumptions  $(\mathcal{A}_1)$ to  $(\mathcal{A}_3)$ hold. Then,  for any $x\in\mathcal{X}$ and any  $j\in\{1,\cdots,q\}$, we have:
\begin{eqnarray*}
\textrm{(i) }\,\,\textrm{IF}(x;\rho_{\gamma\cdot j} ,\mathbb{P}_X )& =&- \frac{\sigma_\gamma^{-2}}{2\kappa_0}\textrm{\large \textbf{1}}_{E_\gamma} (x)\sum_{k=1}^{K} \sum_{\underset{\ell\neq k}{\ell=1}}^{K} < \ \beta_{k}^{(j)}  ,  x_k >_k\,< x_\ell - V_{\ell k}x_k \ , \ \beta_{\ell}^{(j)} \ >_\ell.
\end{eqnarray*}
\noindent\textrm{(ii) }\,\,We suppose, in addition, that    $\rho_1>\rho_2>\cdots >\rho_q$.  Then :
\begin{eqnarray*}
\textrm{IF}( x , \alpha^{(j)}_\gamma ,\mathbb{P}_X) 
&=&-\frac{\sigma_\gamma^{-3}}{2\kappa_0}\textrm{\large \textbf{1}}_{E_\gamma} (x)\, \textrm{IF}( x ; \alpha^{(j)} , \mathbb{P}_X)\\
& &+\sigma_\gamma^{-3}\bigg\{\bigg(\frac{1}{2\kappa_0}-\kappa_1-\kappa_2 \Vert V^{-1/2}x\Vert_\mathcal{X}^2 \bigg)\textrm{\large \textbf{1}}_{E_\gamma}(x)-\kappa_4\bigg\}\beta^{(j)},
\end{eqnarray*}
where  $\textrm{IF}( x ; \alpha^{(j)} , \mathbb{P}_X)$ is given  in (\ref{ifalphaj}).
\end{Theo}

\begin{Rmq}

From this theorem, we recover the results of \cite{crouxdehon} which gives the influence function of  MCD estimator of  LCA of two random vectors. Indeed, using the notation of Remark \ref{rmq1} and (\ref{ifrho2}), we deduce from the previous theorem that, when $K=2$, we have 
\[
\textrm{IF}(X;\rho_{j}^2, \mathbb{P}_X )= -\frac{\sigma_\gamma^{-2}}{2\kappa_0}\textrm{\large \textbf{1}}_{E_\gamma} (x)\left(2\rho_j\,u_j\,v_j-\rho_j^2u_j^2-\rho_j^2v_j^2\right),
\]
what is the result obtained in \cite{crouxdehon}.
\end{Rmq}

\section{Asymptotics for RMSLCA}
\noindent

\noindent In this section  we deal with asymptotic expansion for RMSLCA.  We first establish asymptotic normality for $\widetilde{T}_n$ and then  we derive  the asymptotic distribution  of the  canonical coefficients.  

\bigskip

\begin{Theo}\label{loilim}
Under the assumptions  $(\mathcal{A}_1)$ to  $(\mathcal{A}_3)$, $\sqrt{n}(\widetilde{T}_{n} - T)$ converges in distribution,  as $n \rightarrow +\infty$, to a random variable $U_\gamma$  having a  normal distribution in  $\mathcal L( \mathcal{X})$,  with mean 0 and covariance operator equal to that of the random operator
\begin{eqnarray*}
Z_\gamma &=& \sigma_\gamma^{-2} \kappa_3\textrm{\large \textbf{1}}_{E_\gamma} (X)\, \sum_{k=1}^K \sum_{\underset{\ell\neq k}{\ell=1}}^K \bigg\{- \frac{1}{2}\bigg(\tau_k^* \left(X_k\otimes X_k\right) V_{k\ell}\tau_\ell +  \tau_\ell^*V_{\ell k} \left(X_k\otimes X_k\right)\tau_k\bigg)\\
& & +	\tau_k^*\left(X_\ell\otimes X_k\right) \tau_\ell\bigg\}
+ \left(\sigma_\gamma^{-2} - 1\right) w\left(\Vert V^{-1/2}X\Vert _\mathcal{X} \right) \sum_{k=1}^K \sum_{\underset{\ell\neq k}{\ell=1}}^K\tau_k^*V_{k\ell}\tau_\ell. 
\end{eqnarray*}
where $w\,:\,[0,+\infty[\rightarrow\mathbb{R}$ is the function defined by 
\begin{equation}\label{w}
w(t)= \textrm{\large \textbf{1}}_{[0,r(\gamma)]} (t)\bigg(\kappa_1+\kappa_2\,t^2\bigg)+\kappa_4
\end{equation}
and  $\kappa_1$, $\kappa_2$, $\kappa_3$ and $\kappa_4$ are  given in (\ref{const2}).
\end{Theo}

\bigskip

\noindent This theorem permits to obtain  asymptotic distributions for the canonical coefficients. Let  $\left(\rho^\prime_j\right)_{1\leq j\leq s}$  (with $s\in\mathbb{N}^\ast$) be the decreasing  sequence of distinct eigienvalues of $T$, and    $m_j$ the multiplicity of $\rho^\prime_j$. Putting $\eta_j=\sum_{k=0}^{j-1}m_k$  with $m_0:=0$, we clearly have $\rho_i=\rho^\prime_j$  for any $i\in\{\eta_{j-1}+1,\cdots,\eta_j\}$. We denote by $\Pi_j$ the orthogonal projector from $\mathcal{X}$ onto the eigenspace associated with $\rho^\prime_j$, and by  $\Delta $ the continuous map which associates to each self-adjoint operator $A$ the vector $\Delta (A)$ of its eigenvalues in nonincreasing order. For $j\in\{1,\cdots s\}$, we consider the  $m_j$-dimensional vectors 
\[
\hat{\upsilon}_j=\left(
\begin{array}{c}
\rho^\prime_j\\
\vdots\\
\rho^\prime_j
\end{array}
\right)\,\,\,\,\textrm{ and }\,\,\,\,
\hat{\upsilon}^n_j=\left(
\begin{array}{c}
\widehat{\rho}_{\eta_{j-1}+1\cdot n}\\
\vdots\\
\widehat{\rho}_{\eta_{j}\cdot n}
\end{array}
\right)
\]
from what we consider
\[
\Lambda=\left(
\begin{array}{c}
\upsilon_1\\
\vdots\\
\upsilon_s
\end{array}
\right)
\,\,\,\,\textrm{and}\,\,\,\,
\widehat{\Lambda}_n=\left(
\begin{array}{c}
\hat{\upsilon}^n_1\\
\vdots\\
\hat{\upsilon}^n_s
\end{array}
\right)
.
\]
Then letting $\{e_m\}_{1\leq m\leq q}$ be an orthonormal basis of $\mathcal{X}$, we have:

\bigskip

\begin{Theo}\label{loicoeff}
Under the assumptions  $(\mathcal{A}_1)$ to  $(\mathcal{A}_3)$,  $\sqrt{n}\left(\widehat{\Lambda}_n -\Lambda\right)$ converges in distribution, as $n\rightarrow +\infty$, to the  random vector
\begin{eqnarray*}\label{loivalp}
\zeta=\left(
\begin{array}{c}
\Delta (\Pi_1W_\gamma\Pi_1)\\
\vdots\\
\Delta (\Pi_sW_\gamma\Pi_s)
\end{array}
\right),
\end{eqnarray*}
where $W_\gamma$ is  a random variable having a normal distribution in $\mathcal{L}(\mathcal{X})$, with mean $0$ and covariance operator $\Sigma$ given by:
\begin{eqnarray*}
\Sigma&=&\sum_{1\leq m,r,u,t\leq q}\mathbb{E}\left(\mathcal{Y}_{m,r}\,\mathcal{Y}_{u,t}\right)\,\,(e_m\otimes e_r)\widetilde{\otimes}(e_u\otimes e_t)
\end{eqnarray*}
with
\begin{eqnarray*}
\mathcal{Y}_{m,r}
&=&
\sum_{k=1}^K\sum_{\stackrel{\ell=1}{\ell\neq k}}^K\bigg\{\sigma_\gamma^{-2} \kappa_3\textrm{\large \textbf{1}}_{E_\gamma} (X)\bigg[ - \frac{1}{2}\bigg(<\tau_\ell\beta^{(m)},V_{\ell k}X_k>_\ell<\tau_k\beta^{(r)},X_k>_k\\
& &+<\tau_\ell\beta^{(r)},V_{\ell k}X_k>_\ell <\tau_k\beta^{(m)},X_k>_k\bigg)+<\tau_\ell\beta^{(m)},X_\ell>_\ell<\tau_k\beta^{(r)},X_k>_k\bigg]\\
& & -\left(\sigma_\gamma^{-2} - 1\right)\left( w\left(\Vert V^{-1/2}X\Vert _\mathcal{X} \right)-\kappa_1\gamma -\kappa_2\mu-\kappa_4 \right) <\tau_k\beta^{(r)}, V_{k\ell}\tau_\ell \beta^{(m)}>_k\bigg\},
\end{eqnarray*}
 the function  $w$ being  defined  in (\ref{w}), $\mu=\mathbb{E}\bigg(\textrm{\large \textbf{1}}_{E_\gamma} (X) \,\Vert X\Vert_\mathcal{X}^2\bigg)$ and $\widetilde{\otimes}$ being  the tensor product related to the inner product $<A,B>=\textrm{tr}\left(AB^\ast \right)$.
\end{Theo}

\bigskip

\noindent When the eigenvalues of $T$ are  simple, that is $\rho_1>\rho_2>\cdots >\rho_q$, the  preceding theorem has a simpler statement. We have:

\bigskip

\begin{Cor}\label{loicoeff2}
We suppose that  the assumptions  $(\mathcal{A}_1)$ to  $(\mathcal{A}_3)$ hold and that the canonical coefficient satisfy:  $\rho_1>\rho_2>\cdots >\rho_q$. Then,    $\sqrt{n}\left(\widehat{\Lambda}_n -\Lambda\right)$ converges in distribution, as $n\rightarrow +\infty$, to a random variable having a normal distribution in  $\mathbb{R}^p$ with mean $0$ and covariance matrix $\mathcal{M}=\left(\sigma_{ij}\right)_{1\leq i,j\leq p}$ with:
\[
\sigma_{ij}=\sum_{1\leq m,r,u,t\leq q}\beta^{(i)}_m\beta^{(i)}_r\beta^{(j)}_s\beta^{(j)}_t\mathbb{E}\left(\mathcal{Y}_{m,r}\,\mathcal{Y}_{u,t}\right).
\]
\end{Cor}
\noindent The proof of this corollary  is in all respects similar to that of Corollary 3.1 of \cite{nkiet}, it is then omitted.

\section{Robust test for mutual non-correlation}
In this section  we consider the problem of testing for mutual non-correlation
between  $X_1, X_2, ... , X_K$. This is testing for the null hypothesis 
\begin{center}
$\mathcal H_0 : \forall (k, \ell) \in \left\{1, ... , K \right\}^2,$  $k \neq \ell,$  $V_{k\ell} = 0$
\end{center}
against the alternative  
\begin{center}
$\mathcal H_1 : \exists (k, \ell) \in  \left\{1, ... , K\right\}^2$,   $k \neq \ell$,   $V_{k\ell} \neq 0.$
\end{center}
This testing problem  was already considered in \cite{nkiet}; a test statistic which depends on empirical covariance operator was then  proposed and its asymptotic distribution under the null hypothesis was derived. Since the resulting testing  method may be nonrobust because  its depends on an estimator which is itself nonrobust, it could be interesting to propose a new method that depends instead on a robust estimator of the covariance operator of $X$. Here, we introduce  a test statistic constructed similarly to the one of \cite{nkiet}, but with the MCD estimator of the aforementioned covariance operator. It is then defined  as
\begin{eqnarray}
\widetilde{S}_{n} &=& \displaystyle \sum_{k=2}^{K} \sum_{\ell=1}^{k - 1}\textrm{tr}\left( \pi_{k\ell}\left(\widetilde{T}_{n}\right)\pi_{k\ell}\left(\widetilde{T}_{n}\right)^*\right),
\end{eqnarray}
where  $\widetilde{T}_{n}$ is the estimator given in (\ref{tntilde}) and $\pi_{k\ell}$ is the operator defined as
\[
\pi_{k\ell}\,:\,A\in\mathcal{L}(\mathcal{X})\mapsto \tau_k A\tau_\ell^\ast\in\mathcal{L}(\mathcal{X}_\ell,\mathcal{X}_k).
\] 
\subsection{Asymptotic distribution under the null hypothesis}

\noindent Let us consider
\[
\tau=\frac{\sigma_\gamma^{-4} \kappa_3^2}{q(q+1)}\,\mathbb{E}\bigg(\textrm{\large \textbf{1}}_{E_\gamma} (X) \,\Vert X\Vert_\mathcal{X}^4\bigg).
\]
Then, we have:

\bigskip

\begin{Theo}\label{test}
Suppose that  the assumptions  $(\mathcal{A}_1)$ to  $(\mathcal{A}_3)$ hold. Then, under $\mathcal H_0,$ the sequence $\tau^{-1}n\widetilde{S}_{n}$ converges in distribution, as $n \rightarrow +\infty$, to  $\chi_d^2$, where $d=\sum_{k=1}^K\sum_{\ell =1}^{k-1}p_kp_\ell$ with $p_k=\textrm{dim}(\mathcal{X}_k)$.
\end{Theo}

\bigskip

\begin{Rmq}
For performing this test in practice one has to estimate the unknown parameter $\tau$. This can be done by using an estimate $\widehat{\kappa}_3$ of $\kappa_3$  as in (\ref{const2}) by replacing $r(\gamma)$ (resp. $\sigma_\gamma$) by an estimate $\widehat{r}$ (resp. $\widehat{\sigma}$), and by considering
\[
\widehat{\tau}=\,\frac{\widehat{\sigma}^{-4}\widehat{\kappa}_3^2}{q(q+1)n}\sum_{i=1}^n\textrm{\large \textbf{1}}_{\widehat{E}_\gamma} (X^{(i)}) \,\Vert X^{(i)}\Vert_\mathcal{X}^4,
\]
where $\widehat{E}_\gamma=\left\{x \in  \mathcal{X} , <x ,\widetilde{V}_n^{-1}x >_\mathcal{X} \,\leq\widehat{r}^2\right\}$. A consistent estimator of $r(\gamma)$ is defined in \cite{butler}.

\end{Rmq}
\subsection{Second order influence function}
In order to study robustness properties of the proposed test, we have to derive the influence function related to $\widetilde{S}_n$ under the  null hypothesis. The functional $\mathbb S_{\gamma}$  related to this test statistic is defined as
\[
\mathbb{S}_{\gamma}(\mathbb{P})=\sum_{k=2}^{K}\sum_{\ell=1}^{k-1}\textrm{tr}\left(\pi_{k\ell}\left(\mathbb{T}_{1,\gamma}(\mathbb{P})\right)\pi_{k\ell}\left(\mathbb{T}_{1,\gamma}(\mathbb{P})\right)^*\right)
\]
and, putting $ S_\gamma=\mathbb{S}_{\gamma}(\mathbb{P}_X)$ we, therefore, obtain the related  influence function as
\begin{eqnarray*}
\textrm{IF}(x; S_\gamma, \mathbb{P}_X) &=& \frac{\partial \mathbb S_\gamma(\mathbb{P}_{\varepsilon ,x})}{\partial\epsilon}\bigg\vert_{\varepsilon =0}\\
&=&\sum_{k=2}^{K} \sum_{\ell=1}^{k - 1}\textrm{tr}\left( \pi_{k\ell}\left(\textrm{IF}(x;T,\mathbb{P}_X) \right)\pi_{k\ell}\left( T \right)^*\right)\\
& &+ \sum_{k=2}^{K} \sum_{\ell=1}^{k - 1}\textrm{tr}\left( \pi_{k\ell}\left(T\right)\pi_{k\ell}\left(\textrm{IF}(x;T_{\gamma},\mathbb{P}_X)  \right)^*\right),
\end{eqnarray*}
\noindent
where $\mathbb{P}_{\varepsilon ,x}=(1-\varepsilon)\mathbb{P}_X+\varepsilon \delta_x$ with $\varepsilon\in [0,1]$. Since under  $\mathcal H_0$ we have  $T = 0$ , it follows that  $\textrm{IF}(x; S_\gamma, \mathbb{P}_X)=0$ for all   point $x\in\mathcal{X}$.  In such case, it is necessary to derive the  second order   order influence function of the test statistic, defined as
\[
\textrm{IF}^{(2)}(x; S_\gamma, \mathbb{P}_X)=\frac{\partial^2 \mathbb S_\gamma(\mathbb{P}_{\varepsilon ,x})}{\partial\epsilon^2}\bigg\vert_{\varepsilon =0}.
\]
We have:
\bigskip

\begin{Theo}\label{if2}
We suppose that  the assumptions  $(\mathcal{A}_1)$ to  $(\mathcal{A}_3)$ hold. Then, under $\mathcal H_0$,  the 
second order influence function of $S_\gamma$  is given by
\begin{equation}\label{if2}
\textrm{IF}^{(2)}(x; S_\gamma, \mathbb{P}_X)=\frac{\sigma_\gamma^{-4}}{4\kappa_0^2}\textrm{\large \textbf{1}}_{E_\gamma} (x)\sum_{k=2}^{K} \sum_{\ell=1}^{k - 1}\Vert x_k\Vert  _k^2\,\,\Vert x_\ell\Vert  _\ell^2.
\end{equation}
\end{Theo}
\noindent It is easily  seen that this second order influence function is bounded. Indeed, if $x\in E_\gamma$  then $\Vert x\Vert_\mathcal{X}\leq \Vert V^{1/2}\Vert_\infty\,r(\gamma)$. In addition,   since $\Vert x_k\Vert_k\leq\Vert x\Vert_\mathcal{X}$ we deduce from (\ref{if2})  that
\[
\sup_{x\in\mathcal{X}}\left(\textrm{IF}^{(2)}(x; S_\gamma, \mathbb{P}_X)\right)\leq \frac{\sigma_\gamma^{-4}}{4\kappa_0^2}(K-1)^2\Vert V^{1/2}\Vert_\infty^4\,r^4(\gamma).
\]
This shows that by using   the  MCD estimator of $V$ we have obtained a robust test for mutual non-correlation.

\section{Proofs}
\label{proofs}
\subsection{Proof of Theorem \ref{ift}}
Since
\[
f(V)=\sum_{k=1}^K\tau_k^\ast \tau_kV\tau_k^\ast\tau_k=\sum_{k=1}^K\tau_k^\ast V_k\tau_k=\sum_{k=1}^K\tau_k^\ast \tau_k=\mathbb{I},
\]
where $\mathbb{I}$ is the identity operator of $\mathcal{X}$,  we obtain
\begin{eqnarray*}
\mathbb{T}\left(\mathbb{P}_{\varepsilon ,x}\right)-\mathbb{T}\left(\mathbb{P}_{X}\right)
&=& f\left(\mathbb{V}(\mathbb{P}_{\varepsilon ,x} )\right)^{-1/2}g\left(\mathbb{V}(\mathbb{P}_{\varepsilon ,x} )\right)f\left(\mathbb{V}(\mathbb{P}_{\varepsilon ,x} )\right)^{-1/2}-g(V)\\
&=&\bigg(f\left(\mathbb{V}(\mathbb{P}_{\varepsilon ,x} )\right)^{-1/2}-\mathbb{I}\bigg) g\left(\mathbb{V}(\mathbb{P}_{\varepsilon ,x} )\right)f\left(\mathbb{V}(\mathbb{P}_{\varepsilon ,x} )\right)^{-1/2}\\
& &+\bigg(g\left(\mathbb{V}(\mathbb{P}_{\varepsilon ,x} )-V\right)\bigg)f\left(\mathbb{V}(\mathbb{P}_{\varepsilon ,x} )\right)^{-1/2}\\
&  &+g(V)\bigg(f\left(\mathbb{V}(\mathbb{P}_{\varepsilon ,x} )\right)^{-1/2}-\mathbb{I}\bigg) .
\end{eqnarray*}
where     $\mathbb{P}_{\varepsilon ,x}= (1-\epsilon)\mathbb{P}_X + \varepsilon\delta_x$ with  $\varepsilon \in \left[0; 1\right]$. Then, using the equality 
\begin{equation}\label{decomp}
A^{-1/2}-\mathbb{I}=-A^{-1}(A-\mathbb{I})\left(A^{-1/2}+\mathbb{I}\right)^{-1}
\end{equation}
we obtain:
\begin{eqnarray*}
& &\mathbb{T}\left(\mathbb{P}_{\varepsilon ,x}\right)-\mathbb{T}\left(\mathbb{P}_{X}\right)\\
&=&-f\left(\mathbb{V}(\mathbb{P}_{\varepsilon ,x} )\right)^{-1}f\bigg(\mathbb{V}\left(\mathbb{P}_{\varepsilon ,x} \right)-\mathbb{V}\left(\mathbb{P}_{X}\right)\bigg) \bigg(f\left(\mathbb{V}(\mathbb{P}_{\varepsilon ,x} )\right)^{-1/2}+\mathbb{I}\bigg)^{-1}g\left(\mathbb{V}(\mathbb{P}_{\varepsilon ,x} )\right)f\left(\mathbb{V}(\mathbb{P}_{\varepsilon ,x} )\right)^{-1/2}\nonumber\\
& &+g\bigg(\mathbb{V}(\mathbb{P}_{\varepsilon ,x} )-\mathbb{V}\left(\mathbb{P}_{X}\right)\bigg)f\left(\mathbb{V}(\mathbb{P}_{\varepsilon ,x} )\right)^{-1/2}\nonumber\\
&  &-g(V)f\left(\mathbb{V}(\mathbb{P}_{\varepsilon ,x} )\right)^{-1}f\bigg(\mathbb{V}\left(\mathbb{P}_{\varepsilon ,x} \right)-\mathbb{V}\left(\mathbb{P}_{X}\right)\bigg) \bigg(f\left(\mathbb{V}(\mathbb{P}_{\varepsilon ,x} )\right)^{-1/2}+\mathbb{I}\bigg)^{-1}.
\end{eqnarray*}
Then, from
\[
\textrm{IF}\left(x;V,\mathbb{P}_X \right)=\lim_{\varepsilon\rightarrow 0}\frac{\mathbb{V}\left(\mathbb{P}_{\varepsilon ,x} \right)-\mathbb{V}\left(\mathbb{P}_{X}\right)}{\varepsilon}
\]
and the continuity of the maps $\varepsilon\mapsto \mathbb{V}\left(\mathbb{P}_{\varepsilon ,x}\right)$, $A\mapsto A^{-1}$, $A\mapsto A^{-1/2}$, we deduce that
\begin{eqnarray}\label{relif}
\textrm{IF}\left(x;T,\mathbb{P}_X \right)
&=&-f\left(\mathbb{V}(\mathbb{P}_{X} )\right)^{-1}f(\textrm{IF}\left(x;V,\mathbb{P}_X \right) )\bigg(f\left(\mathbb{V}(\mathbb{P}_{X} )\right)^{-1/2}+\mathbb{I}\bigg)^{-1}g\left(\mathbb{V}(\mathbb{P}_{X} )\right)f\left(\mathbb{V}(\mathbb{P}_{X} )\right)^{-1/2}\nonumber\\
& &+g(\textrm{IF}\left(x;V,\mathbb{P}_X \right))f\left(\mathbb{V}(\mathbb{P}_{X} )\right)^{-1/2}\nonumber\\
&  &-g(V)f\left(\mathbb{V}(\mathbb{P}_{X} )\right)^{-1}f(\textrm{IF}\left(x;V,\mathbb{P}_X \right)) \bigg(f\left(\mathbb{V}(\mathbb{P}_{X} )\right)^{-1/2}+\mathbb{I}\bigg)^{-1}\nonumber\\
&=&-\frac{1}{2}f(\textrm{IF}\left(x;V,\mathbb{P}_X \right) )g\left(V\right)-\frac{1}{2}g\left(V\right)f(\textrm{IF}\left(x;V,\mathbb{P}_X \right) )+g(\textrm{IF}\left(x;V,\mathbb{P}_X \right)).
\end{eqnarray}
Since $\textrm{IF}\left(x;V,\mathbb{P}_X \right))=x\otimes x-V$ (see \cite{critchley}), it follows that
\begin{eqnarray*}
& &\textrm{IF}\left(x;T,\mathbb{P}_X \right)\\
&=&-\frac{1}{2}f\left(x\otimes x\right)g\left(V\right)
-\frac{1}{2}g\left(V\right)f\left( x\otimes x\right)
+g\left(x\otimes x\right)\\
&=&\sum_{k=1}^K\sum_{\stackrel{\ell=1}{\ell\neq k}}^K\bigg(-\frac{1}{2}\sum_{j=1}^K\bigg(\tau_j^\ast\tau_j(x\otimes x)\tau_j^\ast\tau_j\tau_k^\ast V_{k\ell}\tau_\ell
+\tau_k^\ast V_{k\ell}\tau_\ell\tau_j^\ast\tau_j(x\otimes x)\tau_j^\ast\tau_j\bigg)
+\tau_k^\ast\tau_k(x\otimes x)\tau_\ell^\ast\tau_\ell\bigg);
\end{eqnarray*} 
from $\tau_n\tau_m^\ast=\delta_{nm}I_n$, where $\delta$ is the susual Kronecker symbol, and from the equality  $A(x\otimes y)B^\ast=(Bx)\otimes (Ay)$  (see \cite{dauxois}), we deduce that
\begin{eqnarray*}
\textrm{IF}\left(x;T,\mathbb{P}_X \right)
&=&\sum_{k=1}^K\sum_{\stackrel{\ell=1}{\ell\neq k}}^K\bigg(-\frac{1}{2}\bigg(\tau_k^\ast(x_k\otimes x_k)V_{k\ell}\tau_\ell
+\tau_k^\ast V_{k\ell}(x_\ell\otimes x_\ell)\tau_\ell\bigg)
+\tau_k^\ast(x_\ell\otimes x_k)\tau_\ell\bigg)\\
&=&\sum_{k=1}^K\sum_{\stackrel{\ell=1}{\ell\neq k}}^K\bigg(-\frac{1}{2}\bigg(\tau_k^\ast(x_k\otimes x_k)V_{k\ell}\tau_\ell
+\tau_\ell^\ast V_{\ell k}(x_k\otimes x_k)\tau_k\bigg)
+\tau_k^\ast(x_\ell\otimes x_k)\tau_\ell\bigg).
\end{eqnarray*}

\subsection{Proof of Proposition \ref{nonbornee}}
Clearly, $\textrm{IF}(x;T,\mathbb{P}_X)=\textrm{IF}(x;T,\mathbb{P}_X)^\ast$. Then, putting
\[
\theta_{k\ell}=-\frac{1}{2}\tau_k^\ast(x_k\otimes x_k)V_{k\ell}\tau_\ell
-\frac{1}{2}\tau_\ell^\ast V_{\ell k}(x_k\otimes x_k)\tau_k\
+\tau_k^\ast(x_\ell\otimes x_k)\tau_\ell,
\]
we have:
\begin{eqnarray*}
\Vert\textrm{IF}(x;T,\mathbb{P}_X)\Vert_{\mathcal{L}(\mathcal{X})}^2&=& \textrm{tr}\left(\textrm{IF}(x;T,\mathbb{P}_X)^2\right)
= \textrm{tr}\left(\left[\sum_{k=1}^K\sum_{\underset{\ell\neq k}{\ell =1}}^{K}\theta_{k\ell}\right]^2\right)\\
&=& \sum_{k=1}^K \sum_{\underset{\ell\neq k}{\ell =1}}^{K} \bigg\{  \textrm{tr}\left(\theta_{k\ell}^2\right) +   \sum_{\underset{m\neq \ell,\,m\neq k}{m=1}}^{K} \textrm{tr}\left(\theta_{k\ell}\theta_{km}\right) +  \textrm{tr}\left( \theta_{k\ell}\theta_{\ell k}\right) \\
& &+ \sum_{\underset{m\neq k,\,m\neq \ell}{m=1}}^{K} \textrm{tr}\left( \theta_{k\ell}\theta_{\ell m}\right) +  \sum_{\underset{j\neq \ell,\,j\neq k}{j=1}}^{K}\textrm{tr}\left(\theta_{k\ell}\theta_{j\ell}\right)+ \sum_{\underset{j\neq k,\,j\neq\ell}{j=1}}^{K}\textrm{tr}\left(\theta_{k\ell}\theta_{jk}\right)\\
& &+\sum_{\underset{j\neq k,\,j\neq \ell}{j=1}}^{K}\,\,\,  \sum_{\underset{m\neq k,\,m\neq\ell,\,m\neq j}{m=1}}^{K} \textrm{tr}\left(\theta_{k\ell}\theta_{jm}\right)\bigg\}.
\end{eqnarray*}
Using the properties $(a\otimes b)(c\otimes d)=<a,d>c\otimes b$, $A(y\otimes z)B^\ast=(By)\otimes (Az)$, tr$(y\otimes z)=<y,z>$ (see \cite{dauxois}), together with tr$(AB)=\textrm{tr}(BA)$ and $\tau_k\tau_l^\ast=\delta_{kl}I_k$, we obtain:
\begin{eqnarray*}
\textrm{tr}\left(\theta_{k\ell}\theta_{jm}\right)&=& 
 \frac{1}{2}\delta_{k m}\delta_{\ell j}<V_{\ell k}x_k,x_\ell>^2_\ell+\frac{1}{2}\delta_{k j}\delta_{\ell m}\Vert V_{\ell k}x_k\Vert_\ell^2\Vert x_k\Vert^2_k\\
& &-\delta_{k j}\delta_{\ell m}<V_{\ell k}x_k,x_\ell>_\ell\Vert x_k\Vert^2_k -\frac{1}{2}\delta_{k m}\delta_{\ell j}<V_{\ell k}x_k,x_\ell>_\ell\Vert x_k\Vert^2_k\\
& &-\frac{1}{2}\delta_{k m}\delta_{\ell j}<V_{\ell k}x_k,x_\ell>_\ell\Vert x_\ell\Vert^2_\ell+\Vert x_k\Vert^2_k\,\Vert x_\ell\Vert^2_\ell.
\end{eqnarray*}
Clearly,  $\textrm{tr}\left(\theta_{k\ell}\theta_{jm}\right)=0$ if $j\notin\{k,\ell\}$ or if $m\notin\{k,\ell\}$. Hence
\begin{eqnarray}\label{normeif}
\Vert\textrm{IF}(x;T,\mathbb{P}_X)\Vert_{\mathcal{L}(\mathcal{X})}^2
&=&\sum_{k=1}^K \sum_{\underset{\ell\neq k}{\ell =1}}^{K}\bigg\{\frac{1}{2}\Vert V_{\ell k}x_k\Vert_\ell^2\,\Vert x_k\Vert^2_k
-2<V_{\ell k}x_k,x_\ell>_\ell\Vert x_k\Vert^2_k\nonumber\\
& &+\frac{1}{2}<V_{\ell k}x_k,x_\ell>^2_\ell+\Vert x_k\Vert^2_k\,\Vert x_\ell\Vert^2_\ell\bigg\}.
\end{eqnarray}
First, if the $X_k$'s are mutually non-correlated, that is $V_{k\ell}=0$ for any $(k,\ell)\in\{1,\cdots,K\}^2$ with $k\neq\ell$. Then, from (\ref{normeif}) we obtain for any $t\in\mathbb{R}$ and any $x_0=(x^0_1,\cdots,x^0_K)\in\mathcal{X}$ such that $x^0_k\neq 0$ for $k\in\{1,\cdots,K\}$:
\[
\Vert\textrm{IF}(tx_0;T,\mathbb{P}_X)\Vert_{\mathcal{L}(\mathcal{X})}=t^2\sqrt{\sum_{k=1}^K \sum_{\underset{\ell\neq k}{\ell =1}}^{K}\Vert x^0_k\Vert^2_k\,\Vert x^0_\ell\Vert^2_\ell}
\]
from what we deduce that $\lim_{t\rightarrow +\infty}\Vert\textrm{IF}(tx_0;T,\mathbb{P}_X)\Vert_{\mathcal{L}(\mathcal{X})}=+\infty$. Secondly, if    the $X_k$'s are not  mutually non-correlated, there exists a pair $(k_0,\ell_0)\in\{1,\cdots,K\}^2$ such that  $k_0\neq\ell_0$ and $V_{\ell_0 k_0}\neq 0$. Then, considering a vector  $a\in\mathcal{X}_{k_0}-\{0\}$  such that  $V_{\ell_0 k_0}a\neq 0$, we put $x_0=\tau^\ast_{k_0}a$. Clearly, $x_0=(x^0_1,\cdots,x^0_K)$ with $x^0_{k_0}=a$ and  $x^0_{k}=0$ for any $k\in\{1,\cdots,K\}$ such that $k\neq k_0$. Hence, we have
$<V_{\ell k}x^0_k,x^0_\ell>_\ell=0$ and $\Vert x^0_k\Vert^2_k\,\Vert x^0_\ell\Vert^2_\ell=0$ for any pair $(k,\ell)\in\{1,\cdots,K\}^2$ with $k\neq\ell$. Then, from (\ref{normeif}), we deduce that for any $t\in\mathbb{R}$,
\begin{eqnarray*}
\Vert\textrm{IF}(tx_0;T,\mathbb{P}_X)\Vert_{\mathcal{L}(\mathcal{X})}^2
&=&\frac{t^4}{2}\sum_{k=1}^K \sum_{\underset{\ell\neq k}{\ell =1}}^{K}\Vert V_{\ell k}x^0_k\Vert_\ell^2\,\Vert x^0_k\Vert^2_k
=\frac{t^4}{2}\sum_{\underset{\ell\neq k_0}{\ell =1}}^{K}\Vert V_{\ell k_0}x^0_{k_0}\Vert_\ell^2\,\Vert x^0_{k_0}\Vert^2_{k_0}\\
&=&\frac{t^4\Vert a\Vert^2_{k_0}}{2}\sum_{\underset{\ell\neq k_0}{\ell =1}}^{K}\Vert V_{\ell k_0}x^0_{k_0}\Vert_\ell^2\\
&\geq &\frac{t^4\Vert a\Vert^2_{k_0}}{2}\Vert V_{\ell_0 k_0}a\Vert_{\ell_0}^2.
\end{eqnarray*}
Consequently,  $\lim_{t\rightarrow +\infty}\Vert\textrm{IF}(tx_0;T,\mathbb{P}_X)\Vert_{\mathcal{L}(\mathcal{X})}=+\infty$.
\subsection{Proof of Therorem \ref{ifcano1}}
\noindent
\textit{(i)}.  Since
$< \beta^{(j)},\beta^{(k)}>_\mathcal{X} = \delta_{jk}$  for all\  $(j, k ) \in \left\{1, ..., q\right\}^2$, we obtain by applying Lemma 3 of \cite{crouxdehon}:
\begin{eqnarray*}
\textrm{IF}(x;\rho_{j},\mathbb{P}_X ) &=&< \beta^{(j)} ,\textrm{ IF}(x;T,\mathbb{P}_X )\beta^{(j)}>_\mathcal{X}\\
&=& \sum_{k=1}^{K} \sum_{\underset{\ell\neq k}{\ell=1}}^{K} \bigg(-\frac{1}{2} < \beta^{(j)} , \tau_k^* \left(x_k\otimes x_k\right) V_{k\ell}\tau_\ell \beta^{(j)}>_\mathcal{X}\\
& &  -\frac{1}{2}< \beta^{(j)} ,  \tau_\ell^*V_{\ell k} \left(x_k\otimes x_k\right)\tau_k \beta^{(j)}>_\mathcal{X}\\
& &	+ < \beta^{(j)}  , \tau_k^*\left( x_\ell\otimes x_k\right)\tau_\ell\beta^{(j)}>_\mathcal{X}\bigg)\\
&=& \sum_{k=1}^{K} \sum_{\underset{\ell\neq k}{\ell=1}}^{K} \bigg(-\frac{1}{2} <\tau_k \beta^{(j)} , \left(x_k\otimes x_k\right) V_{k\ell}\tau_\ell \beta^{(j)}>_k\\
& &  -\frac{1}{2}<\tau_\ell\beta^{(j)} ,  V_{\ell k} \left(x_k\otimes x_k\right)\tau_k \beta^{(j)}>_\ell\\
& &	+ < \tau_k\beta^{(j)}  , \left( x_\ell\otimes x_k\right)\tau_\ell\beta^{(j)}>_k\bigg)\\
&=&
\sum_{k=1}^{K} \sum_{\underset{\ell\neq k}{l=1}}^{K} \bigg(< \beta_k^{(j)} ,  x_k>_k\,< x_\ell , \beta_\ell^{(j)}>_\ell -< \beta_k^{(j)} ,  x_k >_k\,< x_k , V_{ k\ell} \beta_\ell^{(j)}>_k\bigg)\\
&=&\sum_{k=1}^{K} \sum_{\underset{\ell\neq k}{l=1}}^{K} < \beta_k^{(j)} ,  x_k>_k\,< x_\ell-V_{\ell k}x_k , \beta_\ell^{(j)}>_\ell.
\end{eqnarray*} 
\noindent 
\textit{(ii)}.  Since $f(V(\mathbb{P}_X))=f(V)=\mathbb{I}$, we obtain   by applying  the second part of Lemma  3 in \cite{crouxdehon}:
\begin{eqnarray}\label{ifbetaj}
\textrm{IF}(x;\beta^{(j)},\mathbb{P}_X )
 &=&\sum_{\underset{m\neq j}{m=1}}^{q} \frac{1}{\rho_j - \rho_m}< \beta^{(m)} ,  \textrm{IF}(x;T,\mathbb{P}_X )\beta^{(j)}>_\mathcal{X}  \beta^{(m)} \nonumber \\
& & - \frac{1}{2} < \beta^{(j)} , \textrm{IF}(x;f(V),\mathbb{P}_X )\beta^{(j)}>_\mathcal{X}\beta^{(j)}.
\end{eqnarray}
From the equalities
\[
\textrm{IF}(x;f(V),\mathbb{P}_X )=f(\textrm{IF}(x;V,\mathbb{P}_X ))=f(x\otimes x-V)=f(x\otimes x)-\mathbb{I}=\sum_{k=1}^K\tau_k^\ast (x_k\otimes x_k)\tau_k-\mathbb{I}
\]
it follows
\begin{eqnarray}\label{eq1}
< \beta^{(j)} , \textrm{IF}(x;f(V),\mathbb{P}_X )\beta^{(j)}>_\mathcal{X}&=&\sum_{k=1}^K<\beta^{(j)},\tau_k^\ast (x_k\otimes x_k)\tau_k\beta^{(j)}>_\mathcal{X}-\Vert\beta^{(j)}\Vert_\mathcal{X}^2\nonumber\\
&=&\sum_{k=1}^K<\beta^{(j)}_k, (x_k\otimes x_k)\beta^{(j)}_k>_k-1\nonumber\\
&=&\sum_{k=1}^K<\beta^{(j)}_k, x_k>_k^2-1.
\end{eqnarray}
On the other hand,  similar calculations than in \textit{(ii)} give
\begin{eqnarray}\label{eq2}
& &< \beta^{(m)} ,\textrm{ IF}(x;T,\mathbb{P}_X )\beta^{(j)}>_\mathcal{X}\nonumber\\
&=&
\sum_{k=1}^{K} \sum_{\underset{\ell\neq k}{\ell=1}}^{K}\bigg(< \beta_k^{(m)} , x_k >_k< x_\ell , \beta_\ell^{(j)}>_\ell -\frac{1}{2}< \beta_k^{(m)} , x_k>_k < x_k , V_{k\ell} \beta_\ell^{(j)}>_k\nonumber\\ 
& &-\frac{1}{2}< x_k ,  V_{k\ell}\beta_\ell^{(m)} >_k < x_k , \beta_k^{(j)}  > _k\bigg).
\end{eqnarray}
Introducing (\ref{eq1}) and (\ref{eq2}) in (\ref{ifbetaj}), we obtain
\begin{eqnarray*}
\textrm{IF}(x;\beta^{(j)},\mathbb{P} )& =&
\sum_{k=1}^{K} \sum_{\underset{\ell\neq k}{\ell=1}}^{K}\sum_{\underset{m\neq j}{m=1}}^{q} \frac{1}{\rho_j - \rho_m}\bigg(< \beta_k^{(m)} , x_k >_k< x_\ell , \beta_\ell^{(j)}>_\ell\\
& & -\frac{1}{2}< \beta_k^{(m)} , x_k>_k < x_k , V_{k\ell} \beta_\ell^{(j)}>_k\\ 
& &-\frac{1}{2}< x_k ,  V_{k\ell}\beta_\ell^{(m)} >_k < x_l , \beta_k^{(j)}  > _k\bigg)\beta^{(m)}\\
& &-\frac{1}{2}\left( \sum_{k=1}^K<\beta^{(j)}_k, x_k>_k^2-1\right)\beta^{(j)}.
\end{eqnarray*}
Putting  $\alpha^{(j)}( \mathbb{P}_X) = f(V( \mathbb{P}_X))^{-1/2}\beta^{(j)}( \mathbb{P}_X)$, we have:
\begin{eqnarray*}
\alpha^{(j)}(\mathbb{P}_{\epsilon , x}) - \alpha^{(j)}(\mathbb{P}_X) 
&=& f(V(\mathbb{P}_{\epsilon , x}))^{-1/2}\beta^{(j)}( \mathbb{P}_{\epsilon , x}) - \beta^{(j)}(\mathbb{P}_X) \\
&=& f(V(\mathbb{P}_{\epsilon , x}))^{-1/2}\left(\beta^{(j)}(\mathbb{P}_{\epsilon , x}) - \beta^{(j)}(\mathbb{P}_X)\right)  \\
&  &+  \left(f(V(\mathbb{P}_{\epsilon , x}))^{-1/2}  - \mathbb{I} \right)\beta^{(j)}(\mathbb{P}_X)\\
&=&f(V(\mathbb{P}_{\epsilon , x}))^{-1/2}\left(\beta^{(j)}(\mathbb{P}_{\epsilon , x}) - \beta^{(j)}(\mathbb{P}_X)\right) \\
& &-f(V(\mathbb{P}_{\epsilon , x}))^{-1}\left(f(V(\mathbb{P}_{\epsilon , x})  - V(\mathbb{P}_X)) \right)\left( f(V(\mathbb{P}_{\epsilon , x}))^{-1/2} + \mathbb{I}  \right)^{-1}\beta^{(j)}(\mathbb{P}_X).
\end{eqnarray*}

Therefore  
\begin{eqnarray*}
\textrm{IF}( x , \alpha^{(j)} ,\mathbb{P}_X) &=&\textrm{IF}( x , \beta^{(j)} ,\mathbb{P}_X)  - \  \frac{1}{2} f\left(\textrm{IF}( x , V ,\mathbb{P}_X)\right)\beta^{(j)}\\
&=&\sum_{k=1}^{K} \sum_{\underset{l\neq k}{\ell=1}}^{K}\sum_{\underset{m\neq j}{m=1}}^{q} \frac{1}{\rho_j - \rho_m}\bigg(<\beta_k^{(m)} ,  x_k  >_k< x_\ell , \beta_\ell^{(j)}  >_\ell \\ 
& &-\frac{1}{2}  < \beta_k^{(m)} ,  x_k >_k < x_k , V_{k\ell} \beta_\ell^{(j)}  >_k\\  
& &- \frac{1}{2} < x_k ,  V_{k\ell }\beta_\ell^{(m)}  >_\ell < x_k , \beta_k^{(j)} >_k \bigg)\beta^{(m)}\\
& &-\frac{1}{2}\left( \sum_{k=1}^{K}\left[ \tau_k^*\left(x_k\otimes x_k\right)\tau_k + < \beta_k^{(j)}, x_k>_k^2\mathbb{I}\right]  -  2\mathbb{I}\right)\beta^{(j)}.
\end{eqnarray*}

\subsection{Proof of Theorem \ref{iftmcd}}
It is shown in \cite{crouxhaesbroecka} that under spherical distribution $\mathbb{P}_X^0$, one has 
 \[
\textrm{IF}(x;V_\gamma,\mathbb{P}_X^0)=-(2\kappa_0)^{-1}\textrm{\large \textbf{1}}_{\{\Vert x\Vert_\mathcal{X}^2\leq r(\gamma)\}}\left(x\otimes x\right)+w(\Vert x\Vert_\mathcal{X})\,\,\mathbb{I},
\]
where $w$ is the function defined  in (\ref{w}).
Then  affine equivariant property implies that under elliptical model given in assumtion $(\mathcal{A}_3)$ we have:
\begin{eqnarray}\label{ifv}
\textrm{IF}( x ; V_{\gamma} , \mathbb{P}_X)&=& V^{1/2}\left(\frac{-1}{2\kappa_0}\textrm{\large \textbf{1}}_{E_\gamma}(x)\,(V^{-1/2}x)\otimes (V^{-1/2}x) + w(\Vert (V^{-1/2}x)\Vert_\mathcal{X})\,\,\mathbb{I} \right) V^{1/2}\nonumber\\
&=&\frac{-1}{2\kappa_0}\textrm{\large \textbf{1}}_{E_\gamma}(x)\,\,\,x\otimes x + w(\Vert V^{-1/2}x\Vert_\mathcal{X})\,V.
\end{eqnarray}
Putting $V_\gamma=\mathbb{V}_{1,\gamma}\left(\mathbb{P}_X\right)=\sigma_\gamma^2\,V$, we have   $f\left(\mathbb{V}_{1,\gamma}\left(\mathbb{P}_X\right) \right)= \sigma_\gamma^2\,f(V)=\sigma_\gamma^2\, \mathbb{I}$. Thus
\begin{eqnarray*}
\mathbb{T}_{\gamma}\left(\mathbb{P}_{\varepsilon , x}\right) - \mathbb{T}_{\gamma}\left(\mathbb{P}_X\right)&=&\sigma_\gamma^{-2}\bigg\{A_{\varepsilon,\gamma}^{-1/2}g\left(\mathbb{V}_\gamma\left(\mathbb{P}_{\varepsilon , x}\right)\right)A_{\varepsilon,\gamma}^{-1/2} - g\left(V_\gamma\right)\bigg\}\\
&=&\sigma_\gamma^{-2}\bigg\{\bigg(A_{\varepsilon,\gamma}^{-1/2} - \mathbb{I} \bigg)g\left(\mathbb{V}_\gamma\left(\mathbb{P}_{\varepsilon , x}\right)\right)A_{\epsilon,\gamma}^{-1/2}\\
& &+g\bigg(\mathbb{V}_\gamma\left(\mathbb{P}_{\varepsilon , x}\right) - \mathbb{V}_\gamma\left(\mathbb{P}_X\right)\bigg)A_{\varepsilon ,\gamma}^{-1/2}\\
& &+g\left(V_\gamma\right)\left( A_{\varepsilon ,\gamma}^{-1/2} - \mathbb{I}\right)\bigg\}
\end{eqnarray*}
where $A_{\varepsilon ,\gamma} = \sigma_\gamma^{-2}f\left(\mathbb{V}_\gamma\left(\mathbb{P}_{\varepsilon , x}\right)\right) $. Then using (\ref{decomp}), we obtain
\begin{eqnarray*}
\mathbb{T}_{\gamma}\left(\mathbb{P}_{\varepsilon , x}\right) - \mathbb{T}_{\gamma}\left(\mathbb{P}_X\right)
&=&\sigma_\gamma^{-2}\bigg\{-\sigma_\gamma^{-2}A_{\varepsilon ,\gamma}^{-1}f\bigg( \mathbb{V}_\gamma\left(\mathbb{P}_{\varepsilon , x}\right)  - \sigma_\gamma^2V\bigg)\bigg(A_{\varepsilon ,\gamma}^{-1/2}+ \mathbb{I} \bigg)^{-1}g\left(V_\gamma \right)A_{\varepsilon ,\gamma}^{-1/2}\\
& &+ g\bigg(\mathbb{V}_\gamma\left(\mathbb{P}_{\varepsilon , x}\right)  - \sigma_\gamma^2V \bigg)A_{\varepsilon ,\gamma}^{-1/2}\\
& &- \sigma_\gamma^{-2}g\left(V_\gamma\right)A_{\varepsilon ,\gamma}^{-1}f\bigg( \mathbb{V}_\gamma\left(\mathbb{P}_{\varepsilon , x}\right)  - \sigma_\gamma^2V\bigg)\bigg(A_{\varepsilon ,\gamma}^{-1/2} + \mathbb{I} \bigg)^{-1}\bigg\}.
\end{eqnarray*}
Then, since $\lim_{\varepsilon \rightarrow 0} A_{\varepsilon ,\gamma} 
= \mathbb{I}$, we obtain by using the continuity of the maps   $A\mapsto A^{-1}$ and    $A\mapsto A^{-1/2}$: 
\begin{eqnarray}\label{iftg}
\textrm{IF}( x ; T_\gamma, \mathbb{P}_X) &=&\sigma_\gamma^{-2}\bigg\{-\sigma_\gamma^{-2}\frac{1}{2}f\bigg(\textrm{IF}(x;V_\gamma, \mathbb{P}_X)\bigg) g\left(V_\gamma\right)\nonumber\\
&  &+g\bigg(\textrm{IF}(x;V_\gamma, \mathbb{P}_X)\bigg)-\sigma_\gamma^{-1}\frac{1}{2}g\left(V_\gamma\right)f\bigg(\textrm{IF}(x;V_\gamma, \mathbb{P}_X)\bigg)\bigg\}\nonumber\\\nonumber\\
& =&\sigma_\gamma^{-2}\bigg\{-\displaystyle\frac{1}{2}f\bigg(\textrm{IF}(x;V_\gamma, \mathbb{P}_X)\bigg) g\left(V\right)\nonumber\\
& &+g\bigg(\textrm{IF}(x;V_\gamma, \mathbb{P}_X)\bigg)-\frac{1}{2}\,g\left(V\right)f\bigg(\textrm{IF}(x;V_\gamma, \mathbb{P}_X)\bigg)\bigg\}.
\end{eqnarray}
\noindent
Inserting (\ref{ifv}) in (\ref{iftg}) gives the equality
\begin{eqnarray*}
\textrm{IF}( x ; T_\gamma, \mathbb{P}_X)&=&\sigma_\gamma^{-2}\bigg\{\frac{1}{4\kappa_0}\textrm{\large \textbf{1}}_{E_\gamma}(x)\,f(x\otimes x)g\left(V\right) - \frac{1}{2}w(\Vert V^{-1/2}x\Vert)\,g\left(V\right)\nonumber\\
& &+\frac{1}{4\kappa_0}\textrm{\large \textbf{1}}_{E_\gamma}(x)\,g\left(V\right)f(x\otimes x)     - \frac{1}{2}w(\Vert V^{-1/2}x\Vert)\,g\left(V\right)\nonumber\\
& &-\frac{1}{2\kappa_0}\textrm{\large \textbf{1}}_{E_\gamma}(x)\,g(x\otimes x) + w(\Vert V^{-1/2}x\Vert)\,g\left(V\right)\bigg\}\\
&=&\frac{\sigma_\gamma^{-2}}{2\kappa_0}\textrm{\large \textbf{1}}_{E_\gamma}(x)\,\bigg\{\frac{1}{2}f(x\otimes x)\,g\left(V\right)+\frac{1}{2}g\left(V\right)f(x\otimes x)
-g(x\otimes x)\bigg\}\\
&=&-\frac{\sigma_\gamma^{-2}}{2\kappa_0}\textrm{\large \textbf{1}}_{E_\gamma}(x)\,\textrm{IF}( x ; T , \mathbb{P}_X).
\end{eqnarray*}
\subsection{Proof of Proposition \ref{iftmcdbornee}}
If $x\in\E_\gamma$, then $\Vert x\Vert_\mathcal{X}\leq \Vert V^{1/2}\Vert_\infty\,\Vert  V^{-1/2}x\Vert_\mathcal{X}\leq \Vert V^{1/2}\Vert_\infty\,r(\gamma)$. On the other hand,
\begin{eqnarray*}
\Vert \textrm{IF}(x;T _{\gamma}, \mathbb{P}_X)\Vert_{\infty} &\leq&\frac{\sigma_\gamma^{-2}}{2\vert\kappa_0\vert}\textrm{\large \textbf{1}}_{E_\gamma}(x)\,\sum_{k=1}^K \sum_{\underset{\ell\neq k}{\ell=1}}^K \bigg\{ \frac{1}{2}\Vert\tau_k^*\Vert_\infty\,\Vert x_k\otimes x_k\Vert _\infty\,\Vert V_{k\ell}\Vert_\infty\,\Vert\tau_\ell\Vert_\infty \\
& &+\frac{1}{2}\Vert \tau_\ell^* \Vert_\infty\,\Vert V_{\ell k}\Vert_\infty\,\Vert x_k\otimes x_k\Vert_\infty\,\Vert\tau_k\Vert_\infty\,+\Vert\tau_k^*\Vert_\infty\,\Vert x_\ell\otimes x_k\Vert_\infty\,\,\Vert\tau_\ell\Vert_\infty\bigg\} .
\end{eqnarray*}
It is easy to check that
\[
\Vert\tau_k^*\Vert_\infty\leq 1,\,\,\,\Vert\tau_\ell\Vert_\infty\leq 1,\,\,\,\Vert x_k\otimes x_k\Vert_\infty\leq \Vert x_k\Vert_k^2\leq \Vert x\Vert_\mathcal{X}^2
,\,\,\,\Vert x_\ell\otimes x_k\Vert_\infty\leq \Vert x_k\Vert_k\Vert x_\ell\Vert_\ell\leq \Vert x\Vert_\mathcal{X}^2
\]
and
\[
\Vert V_{\ell k}\Vert_\infty=\Vert \tau_k V\tau_{\ell}^\ast\Vert_\infty\leq \Vert\tau_k\Vert_\infty\,\Vert\tau_\ell^*\Vert_\infty\, \Vert V\Vert_\infty\leq \Vert V\Vert_\infty,\,\,\,
\Vert V_{k\ell }\Vert_\infty\leq \Vert V\Vert_\infty.
\]
Hence, for any $x\in\mathcal{X}$,   we have
\begin{eqnarray*}
\Vert \textrm{IF}(x;T _{\gamma}, \mathbb{P}_X)\Vert_{\infty} &\leq&\frac{\sigma_\gamma^{-2}}{2\vert\kappa_0\vert}K(K-1)\bigg(\Vert V\Vert_\infty+1\bigg)\textrm{\large \textbf{1}}_{E_\gamma}(x)\,\Vert x\Vert_\mathcal{X}^2\\
&\leq&\frac{\sigma_\gamma^{-2}}{2\vert\kappa_0\vert}K(K-1)\bigg(\Vert V\Vert_\infty+1\bigg)\Vert V^{1/2}\Vert_\infty^2\,r^2(\gamma).
\end{eqnarray*}
\subsection{Proof of Theorem \ref{ifcanomcd}}
\noindent (i). 
From  Lemma 3 in \cite{crouxdehon}  we obtain
\begin{eqnarray*}
\textrm{IF}( x ; \rho_{\gamma .j} , \mathbb{P}_X)&=&< \beta^{(j)} \  , \textrm{IF}( x ; T_\gamma , \mathbb{P}_X)\beta^{(j)} >_{\mathcal X}\\
&=&-\frac{\sigma_\gamma^{-2}}{2\kappa_0}\textrm{\large \textbf{1}}_{E_\gamma}(x)\,< \beta^{(j)},\textrm{ IF}( x ; T , \mathbb{P}_X)\beta^{(j)} >_{\mathcal X}\\
&=& -\frac{\sigma_\gamma^{-2}}{2\kappa_0}\textrm{\large \textbf{1}}_{E_\gamma}(x)\, \sum_{k=1}^{K} \sum_{\underset{\ell\neq k}{\ell=1}}^{K} < \ \beta_{k }^{(j)}  ,  x_k >_k\,< x_\ell - V_{\ell k}x_k \ , \ \beta_{\ell }^{(j)} \ >_\ell.
\end{eqnarray*}
\noindent (ii).  Since $f(\sigma_\gamma^{-2}V_\gamma(\mathbb{P}_X))=f(V)=\mathbb{I}$, we obtain   by applying  the second part of Lemma  3 in \cite{crouxdehon}:
\begin{eqnarray}\label{ifbetamcd}
\textrm{IF}( x ; \beta_\gamma^{(j)} , \mathbb{P}_X)&=&\sum_{\underset{m\neq j}{m=1}}^{q}  \frac{1}{\rho_{j} - \rho_{m}}< \beta^{(m)},\textrm{IF}( x ; T_\gamma , \mathbb{P}_X)\beta^{(j)} >_{\mathcal X}\beta^{(m)}\nonumber\\
& & -\frac{1}{2} < \beta^{(j)} ,\textrm{IF}( x ; f\left(\sigma_\gamma^{-2}V_\gamma\right) , \mathbb{P}_X)\beta^{(j)} >_{\mathcal X}\beta^{(j)}
\end{eqnarray}
\noindent
Further, $\textrm{IF}( x ; f\left(\sigma_\gamma^{-2}V_\gamma\right) , \mathbb{P}_X)=\sigma_\gamma^{-2}f\left(\textrm{IF}\left( x ; V_\gamma,\left(\mathbb{P}_X\right)\right)\right)$
and from (\ref{ifv}) it follows
\begin{eqnarray*}
\textrm{IF}( x ; f\left(\sigma_\gamma^{-2}V_\gamma\right) , \mathbb{P}_X)&=&-\frac{\sigma_\gamma^{-2}}{2\kappa_0}\textrm{\large \textbf{1}}_{E_\gamma}(x)\,f\left(x\otimes x\right)  + \sigma_\gamma^{-2}w(\Vert V^{-1/2}x\Vert_\mathcal{X}) f\left(V\right)\\
&=&-\frac{\sigma_\gamma^{-2}}{2\kappa_0}\textrm{\large \textbf{1}}_{E_\gamma}(x)\,\bigg(f\left(x\otimes x\right)-\mathbb{I}\bigg)  +\sigma_\gamma^{-2} \bigg(-\frac{1}{2\kappa_0}\textrm{\large \textbf{1}}_{E_\gamma}(x)+w(\Vert V^{-1/2}x\Vert_\mathcal{X}) \bigg)\mathbb{I}\\
&=&-\frac{\sigma_\gamma^{-2}}{2\kappa_0}\textrm{\large \textbf{1}}_{E_\gamma}(x)\,\textrm{IF}( x ; f\left(V\right) , \mathbb{P}_X)  \\
& &+\sigma_\gamma^{-2}\bigg\{ \textrm{\large \textbf{1}}_{E_\gamma}(x)\bigg(-\frac{1}{2\kappa_0}+\kappa_1+\kappa_2 \Vert V^{-1/2}x\Vert_\mathcal{X}^2 \bigg)+\kappa_4\bigg\}\mathbb{I},
\end{eqnarray*}
where $w$ is the function defined in (\ref{w}).
This equality together with  (\ref{ifbetamcd}) and Theorem \ref{iftmcd} imply
\begin{eqnarray*}
\textrm{IF}( x ; \beta_\gamma^{(j)} , \mathbb{P}_X)&=&-\frac{\sigma_\gamma^{-2}}{2\kappa_0}\textrm{\large \textbf{1}}_{E_\gamma} (x)\, \sum_{\underset{m\neq j}{m=1}}^{q}  \frac{1}{\rho_{j} - \rho_{m}}< \beta_{\gamma}^{(m)},\textrm{IF}( x ; T , \mathbb{P}_X)\beta^{(j)} >_{\mathcal X}\beta^{(m)}\nonumber\\
& & +\frac{\sigma_\gamma^{-2}}{4\kappa_0}\textrm{\large \textbf{1}}_{E_\gamma} (x)\, < \beta^{(j)} ,\textrm{IF}( x ; f\left(V\right) , \mathbb{P}_X)\beta^{(j)} >_{\mathcal X}
\beta^{(j)}\\
& & +\sigma_\gamma^{-2}\bigg\{\bigg(\frac{1}{4\kappa_0}-\frac{\kappa_1}{2}-\frac{\kappa_2 }{2}\Vert V^{-1/2}x\Vert_\mathcal{X}^2 \bigg)\textrm{\large \textbf{1}}_{E_\gamma}(x)-\frac{\kappa_4}{2}\bigg\}\beta^{(j)}\\
&=&-\frac{\sigma_\gamma^{-2}}{2\kappa_0}\textrm{\large \textbf{1}}_{E_\gamma} (x)\, \textrm{IF}( x ; \beta^{(j)} , \mathbb{P}_X)\\
& &+\sigma_\gamma^{-2}\bigg\{\bigg(\frac{1}{4\kappa_0}-\frac{\kappa_1}{2}-\frac{\kappa_2 }{2}\Vert V^{-1/2}x\Vert_\mathcal{X}^2 \bigg)\textrm{\large \textbf{1}}_{E_\gamma}(x)-\frac{\kappa_4}{2}\bigg\}\beta^{(j)}.
\end{eqnarray*}
\noindent
On the other hand, since
\[
\alpha^{(j)}_\gamma(\mathbb{P}_X) =\bigg( f(\mathbb{V}_\gamma(\mathbb{P}_X))\bigg)^{-1/2}\beta_\gamma^{(j)}(\mathbb{P}_X) = f(V_\gamma)^{-1/2}\beta_\gamma^{(j)}(\mathbb{P}_X)= \sigma_\gamma^{-1}\beta_\gamma^{(j)}(\mathbb{P}_X),
\] 
it follows
\begin{eqnarray*}
\alpha_\gamma^{(j)}(\mathbb{P}_{\varepsilon , x}) - \alpha_\gamma^{(j)}(\mathbb{P}_X) &=& f(V_\gamma(\mathbb{P}_{\varepsilon , x}))^{-1/2}\beta_\gamma^{(j)}(\mathbb{P}_{\epsilon , x}) - \sigma_\gamma^{-1}\beta_\gamma^{(j)}(\mathbb{P}_X) \nonumber\\
&=& \sigma_\gamma^{-1}\bigg\{A_{ \varepsilon , \gamma}^{-1/2}\beta_\gamma^{(j)}(\mathbb{P}_{\epsilon , x}) - \beta_\gamma^{(j)}(\mathbb{P}_X)\bigg\} \nonumber\\
&=& \sigma_\gamma^{-1}\bigg\{A_{ \varepsilon , \gamma}^{-1/2}\left(\beta_\gamma^{(j)}(\mathbb{P}_{\epsilon , x}) - \beta_\gamma^{(j)}(\mathbb{P}_X)\right)  +  \left(A_{ \varepsilon , \gamma}^{-1/2}  - \mathbb{I} \right)\beta_\gamma^{(j)}(\mathbb{P}_X)\bigg\},
\end{eqnarray*}
where  $A_{\varepsilon ,\gamma} = \sigma_\gamma^{-2}f\left(\mathbb{V}_\gamma\left(\mathbb{P}_{\varepsilon , x}\right)\right) $. Then using (\ref{decomp}), we obtain:
\begin{eqnarray}\label{ifalphamcd}
\alpha_\gamma^{(j)}(\mathbb{P}_{\epsilon , x}) - \alpha_\gamma^{(j)}(\mathbb{P}_X)&=& \sigma_\gamma^{-1}\bigg\{A_{\varepsilon , \gamma}^{-1/2}\left(\beta_\gamma^{(j)}(\mathbb{P}_{\varepsilon , x}) - \beta_\gamma^{(j)}(\mathbb{P}_X)\right)\nonumber\\
&-&A_{ \varepsilon , \gamma}^{-1}\left( A_{ \varepsilon , \gamma}  - \mathbb{I} \right)\left( A_{ \varepsilon , \gamma}^{-1/2} + \mathbb{I}  \right)^{-1}\beta_\gamma^{(j)}(\mathbb{P}_X)\nonumber\\
&=&\sigma_\gamma^{-1}\bigg\{A_{\varepsilon , \gamma}^{-1/2}\left(\beta_\gamma^{(j)}(\mathbb{P}_{\varepsilon  , x}) - \beta^{(j)}(\mathbb{P})\right)\\
&-&\sigma_\gamma^{-2}A_{ \varepsilon , \gamma}^{-1/2}\left(f(V_\gamma(\mathbb{P}_{\varepsilon , x})  - V_\gamma(\mathbb{P}_X)) \right)\left( A_{ \varepsilon , \gamma}^{-1/2} + \mathbb{I}  \right)^{-1}\beta_\gamma^{(j)}(\mathbb{P}_X)\bigg\}\nonumber.
\end{eqnarray}
\noindent
From the continuity of the maps   $A\mapsto A^{-1} , $   $A\mapsto A^{-1/2} $, and the  equality $
\lim_{\varepsilon \rightarrow 0} A_{\varepsilon ,\gamma} = \mathbb{I}$,  we deduce from (\ref{ifalphamcd}) that  
\begin{eqnarray*}
\textrm{IF}( x ;\alpha_\gamma^{(j)} ,\mathbb{P}_X) &=& \lim_{\epsilon\rightarrow 0}\displaystyle\frac{\alpha_\gamma^{(j)}(\mathbb{P}_{\varepsilon , x}) - \alpha_\gamma^{(j)}(\mathbb{P}_X)}{\varepsilon}\\
&=& \frac{1}{\sigma_\gamma}\bigg\{\textrm{IF}( x ; \beta_\gamma^{(j)} ,\mathbb{P}_X)   - \frac{\sigma_\gamma^{-2}}{2} f\left(\textrm{IF}( x ; V_\gamma ,\mathbb{P}_X)\right)\beta^{(j)}\bigg\}\\
&=&-\frac{\sigma_\gamma^{-3}}{2\kappa_0}\textrm{\large \textbf{1}}_{E_\gamma} (x)\, \textrm{IF}( x ; \beta^{(j)} , \mathbb{P}_X)\\
& &+2\sigma_\gamma^{-3}\bigg\{\bigg(\frac{1}{4\kappa_0}-\frac{\kappa_1}{2}-\frac{\kappa_2 }{2}\Vert V^{-1/2}x\Vert_\mathcal{X}^2 \bigg)\textrm{\large \textbf{1}}_{E_\gamma}(x)-\frac{\kappa_4}{2}\bigg\}\beta^{(j)}\\
& &+\frac{\sigma_\gamma^{-3}}{4\kappa_0}\textrm{\large \textbf{1}}_{E_\gamma} (x)\, \textrm{IF}( x ; f(V) , \mathbb{P}_X)\beta^{(j)}\\
& =&-\frac{\sigma_\gamma^{-3}}{2\kappa_0}\textrm{\large \textbf{1}}_{E_\gamma} (x)\, \textrm{IF}( x ; \alpha^{(j)} , \mathbb{P}_X)\\
& &+\sigma_\gamma^{-3}\bigg\{\bigg(\frac{1}{2\kappa_0}-\kappa_1-\kappa_2 \Vert V^{-1/2}x\Vert_\mathcal{X}^2 \bigg)\textrm{\large \textbf{1}}_{E_\gamma}(x)-\kappa_4\bigg\}\beta^{(j)}.
\end{eqnarray*}
\subsection{Proof of Theorem \ref{loilim}}
\subsubsection{A preliminary  lemma}

\noindent The following lemma  gives the asymptotic distribution of the random variable
\begin{equation}\label{hn}
\widehat{H}_{n\gamma}=\sqrt{n}\left(\widetilde{V}_{n}-\sigma_\gamma^2\,V\right).  
\end{equation}

\begin{Lem}\label{loihn}
We assume that assumptions ($\mathcal{A}_1$)  to  ($\mathcal{A}_3$) hold.  Then, $\widehat{H}_{n\gamma} $ converges in distribution  in $\mathcal{L}(\mathcal{X})$, as $n\rightarrow +\infty$, to a random variable $H_\gamma$ having a  normal distribution $N(0,\Lambda)$, where $\Lambda$ is the covariance operator of 
\[
\mathcal{Z}=\kappa_3\, \textrm{\large \textbf{1}}_{E_\gamma} (X)\,  X\otimes  X+w(\Vert  V^{-1/2} X\Vert_\mathcal{X})\,V
\] 
and $w$ is the function given in (\ref{w}).
\end{Lem}
\noindent\textit{Proof}. Using affine equivariant property, we deduce from Eq. (A.25) in \cite{cator} that:
\begin{eqnarray*}
\widehat{H}_{n\gamma}
&=&V^{1/2}\bigg(\frac{1}{\sqrt{n}}\sum_{i=1}^n\bigg(\frac{v (\Vert  V^{-1/2} X^{(i)}\Vert_\mathcal{X})}{\Vert V^{-1/2} X^{(i)}\Vert^2_\mathcal{X}}  ( V^{-1/2} X^{(i)})\otimes ( V^{-1/2} X^{(i)})\\
& &\hspace{3cm}+w(\Vert V^{-1/2} X^{(i)}\Vert_\mathcal{X})\,\mathbb{I}\bigg)+o_P(1)\bigg)V^{1/2}\\
&=&\frac{1}{\sqrt{n}}\sum_{i=1}^n\left(\frac{v (\Vert  V^{-1/2} X^{(i)}\Vert_\mathcal{X})}{\Vert V^{-1/2} X^{(i)}\Vert^2_\mathcal{X}}\,  X^{(i)}\otimes  X^{(i)} +w(\Vert  V^{-1/2} X^{(i)}\Vert_\mathcal{X})\,V\right)+o_P(1)\\
&=&\widehat{W}_n+o_P(1)
\end{eqnarray*}
where $\widehat{W}_n=n^{-1/2}\sum_{i=1}^n\mathcal{Z}_i$, with 
\[
\mathcal{Z}_i=\frac{v (\Vert  V^{-1/2} X^{(i)}\Vert_\mathcal{X})}{\Vert V^{-1/2}X^{(i)}\Vert^2_\mathcal{X}}\,   X^{(i)}\otimes  X^{(i)} +w(\Vert  V^{-1/2} X^{(i)}\Vert_\mathcal{X})\,V,
\]
and  $v\,:\,[0+\infty[\rightarrow [0+\infty[$ is the function defined by  $v(t)=\kappa_3\,\textrm{\large \textbf{1}}_{[0,r(\gamma)]}(t)\,t^2$.  Slustky's theorem permits to conclude that  $\widehat{H}_{n\gamma}$ has the same limiting distribution than  $\widehat{W}_n$, which can be obtained by using central limit  theorem. For doing that, we will first
show that  $\mathbb{E}(\mathcal{Z}_i)=0$. Putting $Y^{(i)}=V^{-1/2} X^{(i)}$, we have 
\begin{equation}\label{ezi}
\mathbb{E}(\mathcal{Z}_i)=V^{1/2}\,\mathbb{E}\bigg(\frac{v(\Vert  Y^{(i)}\Vert_\mathcal{X})}{\Vert Y^{(i)}\Vert^2_\mathcal{X}}\,   Y^{(i)}\otimes  Y^{(i)} +w(\Vert Y^{(i)}\Vert_\mathcal{X})\,\mathbb{I}\bigg)V^{1/2},
\end{equation}
and since   $Y^{(i)}$ has a spherical distribution, we deduce from \cite{cator} (see p. 2387) that $\mathbb{E}\left(v(\Vert Y^{(i)}\Vert_\mathcal{X})+q\,w(\Vert Y^{(i)}\Vert_\mathcal{X})\right)=0$,   $\mathbb{E}\left(v(\Vert Y^{(i)}\Vert_\mathcal{X})^2\right)<+\infty$, and $\mathbb{E}\left(w(\Vert Y^{(i)}\Vert_\mathcal{X})^2\right)<+\infty$. Therefore   $\mathbb{E}\left(w(\Vert Y^{(i)}\Vert_\mathcal{X}\right)=-q^{-1}  \mathbb{E}\left(v(\Vert Y^{(i)}\Vert_\mathcal{X}\right)$, and (\ref{ezi}) becomes:
\begin{eqnarray}\label{ezi2}
\mathbb{E}(\mathcal{Z}_i)&=&V^{1/2}\,\bigg(\mathbb{E}\left(\frac{v (\Vert  Y^{(i)}\Vert_\mathcal{X})}{\Vert Y^{(i)}\Vert^2_\mathcal{X}}\,   Y^{(i)}\otimes  Y^{(i)}\right)-\frac{1}{q}\mathbb{E}(v(\Vert Y^{(i)}\Vert_\mathcal{X}))\,\mathbb{I}\bigg)V^{1/2}\nonumber\\
&=&V^{1/2}\,\bigg(\mathbb{E}\left(\kappa_3\textrm{\large \textbf{1}}_{[0,r(\gamma)]}(\Vert Y^{(i)}\Vert_\mathcal{X})\,   Y^{(i)}\otimes  Y^{(i)}\right)-\frac{1}{q}\mathbb{E}(v(\Vert Y^{(i)}\Vert_\mathcal{X}))\,\mathbb{I}\bigg)V^{1/2}.
\end{eqnarray}
From the proof of Theorem 4.2 in \cite{cator} (see p. 2386) we have
\[
\mathbb{E}\left(\kappa_3\textrm{\large \textbf{1}}_{[0,r(\gamma)]}(\Vert Y^{(i)}\Vert_\mathcal{X})\,   Y^{(i)}\otimes  Y^{(i)}\right)=\frac{1}{q}\mathbb{E}\left(\kappa_3\textrm{\large \textbf{1}}_{[0,r(\gamma)]}(\Vert Y^{(i)}\Vert_\mathcal{X})\,  \Vert Y^{(i)}\Vert_\mathcal{X}^2\right)\,\mathbb{I}
=\frac{1}{q}\mathbb{E}\left(v(\Vert Y^{(i)}\Vert_\mathcal{X})\right)\,\mathbb{I}.
\]
Then,  (\ref{ezi2})  implies  $\mathbb{E}(\mathcal{Z}_i)=0$. Now, using the central limit theorem we conclude that $\widehat{W}_n$ converges in distribution, as $n\rightarrow +\infty$, to a normal distribution $N(0,\Lambda)$ in $\mathcal{L}(\mathcal{X})$, where $\Lambda$ is the covariance operator of 
\begin{eqnarray}\label{rvz}
\mathcal{Z}&=&\frac{v (\Vert  V^{-1/2} X\Vert_\mathcal{X})}{\Vert V^{-1/2}X\Vert^2_\mathcal{X}}\,   X\otimes  X+w(\Vert  V^{-1/2} X\Vert_\mathcal{X})\,V\nonumber\\
&=&\kappa_3\, \textrm{\large \textbf{1}}_{E_\gamma} (X)\,  X\otimes  X+w(\Vert  V^{-1/2} X\Vert_\mathcal{X})\,V.
\end{eqnarray}
\subsubsection{Proof of the theorem}
Arguing as  in  the proof of Theorem 3.2   in \cite{nkiet} (see p. 203), we have the equality   $\sqrt{n}(\widetilde{T}_{n } - T) =\widehat{\varphi}_n\left(\widehat{H}_{n\gamma}\right)$, where
 $\widehat{H}_{n\gamma}$ is given  in (\ref{hn})  and $\widehat{\varphi}_{n}$  is the random operator from  $\mathcal L( \mathcal{H})$ to itself defined  by 
\begin{eqnarray*}
\widehat{\varphi}_{n}(A)&= &-\sigma_\gamma^{-1}f(\widetilde{V}_{n})^{-1}f(A)\left(\sigma_\gamma f(\widetilde{V}_{n})^{-1/2} +  \mathbb{I}\right)^{-1}g(\widetilde{V}_n)f(\widetilde{V}_n)^{-1/2} + \sigma_\gamma^{-1}g(A)f(\widetilde{V}_n)^{-1/2} \\
& &  -g(V)f(\widetilde{V}_n)^{-1}f(A)\left(\sigma_\gamma f(\widetilde{V}_{n})^{-1/2} +  \mathbb{I}\right)^{-1}.
\end{eqnarray*}
Considering the linear map $\varphi_\gamma$ from $\mathcal{X}$ to itself defined as
\[
\varphi_\gamma (A)=\sigma_\gamma^{-2}\left( -\frac{1}{2}f(A)g(V)  + g(A) -\frac{1}{2}g(V)f(A)\right),
\] 
and denoting by  $\Vert\cdot\Vert_{\infty}$ and  $\Vert\cdot\Vert_{\infty \infty}$ the  norms  of  $\mathcal L( \mathcal{X})$\  and \  $\mathcal L(\mathcal L( \mathcal{X}))$,  respectively, defined  by $\Vert A\Vert_{\infty}=\sup_{x \in  \mathcal{X}-\{0\}}\Vert Ax\Vert_ \mathcal{X}/\Vert x\Vert_ \mathcal{X}$ and $\Vert Q\Vert_{\infty \infty}=\sup_{B \in \mathcal L( \mathcal{X})-\{0\}}\Vert Q(B)\Vert_{\infty}/\Vert B\Vert_{\infty}$,  we have :

\begin{eqnarray}\label{ineg1}
 \Vert\widehat{\varphi}_n(\widehat{H}_{n\gamma})-\varphi_\gamma (\widehat{H}_{n\gamma})\Vert_{\infty} 
& \leq&\Vert\widehat{\varphi}_n-\varphi_\gamma \Vert_{\infty\infty}    \, \Vert\widehat{H}_{n\gamma}\Vert_{\infty}
\end{eqnarray}
and 
\begin{eqnarray}\label{ineg2}
\Vert\widehat{\varphi}_n-\varphi_\gamma \Vert_{\infty\infty}  &\leq& \Bigg(\sigma_\gamma^{-1}\Vert f\Vert_{\infty\infty}\,\Vert\left(\sigma_\gamma f(\widetilde{V}_n)^{-1/2} +  \mathbb{I}\right)^{-1}g(\widetilde{V}_{n})f(\widetilde{V}_n)^{-1/2}\Vert_{\infty}\nonumber\\
& &\hspace{1cm}+ \frac{\sigma_\gamma^{-1}}{2}\Vert f\Vert _{\infty \infty} \Vert g(V)\Vert_{\infty} +\sigma_\gamma^{-1} \Vert g\Vert _{\infty \infty}\nonumber\\
& & \hspace{1cm}+ \Vert f\Vert _{\infty \infty}\Vert g(V)\Vert_{\infty}\Vert \left(\sigma_\gamma f(\widetilde{V}_n)^{-1/2} +\mathbb{I}\right)^{-1}\Vert_{\infty} \bigg)\Vert f(\widetilde{V}_n)^{-1/2} -\sigma_\gamma^{-1}\mathbb{I}\Vert_{\infty} \nonumber \\
& +&\bigg( \sigma_\gamma^{-3}\Vert f\Vert _{\infty \infty} \Vert g(\widetilde{V}_{n})f(\widetilde{V}_n)^{-1/2}\Vert_{\infty}\nonumber\\
&  &\hspace{1cm}+ \Vert f\Vert _{\infty \infty}\Vert g(V)\Vert_{\infty}\bigg)\Vert \sigma_\gamma f(\widetilde{V}_n)^{-1/2} +\mathbb{I})^{-1}-\frac{1}{2}\mathbb{I}\Vert_{\infty}\nonumber \\
& +&  \frac{\sigma_\gamma^{-3}}{2} \Vert f\Vert _{\infty\infty} \Vert g\Vert _{\infty \infty}  \Vert f(\widetilde{V}_n)^{-1/2}\Vert_{\infty}\Vert\widetilde{V}_{n}-\sigma_\gamma^2\,V\Vert_{\infty} .
\end{eqnarray}
Lemma \ref{loihn} implies that  $\widetilde{V}_{n}$ converges in probability to $\sigma_\gamma^2\,V$, as $n\rightarrow +\infty$. Then,  using the  continuity of maps $f$,  $g$,  $A\mapsto A^{-1}$ and $A\mapsto A^{-1/2}$  we deduce that $f(\widetilde{V}_n)$ (resp.  $f(\widetilde{V}_n)^{-1}$; resp.  $f(\widetilde{V}_n)^{-1/2}$; resp.  $g(\widetilde{V}_n)$)  converges  in probability, as $n\rightarrow +\infty$,  to  $\sigma_\gamma^2\mathbb{I}$ (resp.  $\sigma_\gamma^{-2}\mathbb{I}$; resp.  $\sigma_\gamma^{-1}\mathbb{I}$;
resp.     $\sigma_\gamma^2\,g(V)$). Consequently, from (\ref{ineg1})  and (\ref{ineg2})  we deduce that  $\widehat{\varphi}_n(\widehat{H}_{n\gamma})-\varphi_\gamma(\widehat{H}_{n\gamma})$ converge in probability to $0$ as $n\rightarrow +\infty$. Slutsky's theorem allows to conclude that $\widehat{\varphi}_n(\widehat{H}_{n\gamma})$ and $\varphi_\gamma(\widehat{H}_{n\gamma})$ both converge to the  same   distribution, that is the distribution of   $\varphi_\gamma (H_\gamma)$.    Since $\varphi_\gamma$ is linear this distribution is the normal distribution with mean equal to $0$ and covariance operator equal to that of the random variable:   
\begin{eqnarray*}
Z_\gamma=\varphi_\gamma\left(\mathcal{Z}\right) &=& \kappa_3\, \textrm{\large \textbf{1}}_{E_\gamma} (X)\,  \varphi_\gamma(X\otimes  X)+w(\Vert  V^{-1/2} X\Vert_\mathcal{X})\,\varphi_\gamma(V).
\end{eqnarray*}
Besides
\begin{eqnarray*}
\varphi_\gamma\left(X\otimes X\right)&=& \sigma_\gamma^{-2}\sum_{k=1}^K \sum_{\underset{\ell\neq k}{\ell=1}}^K - \frac{1}{2}\left(\tau_k^* \left(X_k\otimes X_k\right) V_{k\ell}\tau_\ell  + \tau_\ell^*V_{\ell k} \left(X_k\otimes X_k\right)\tau_k\right) +	\tau_k^*\left(X_\ell\otimes X_k\right) \tau_\ell,
\end{eqnarray*}
and from   $f\left(V\right)= \mathbb{I}$,  it follows:
\[
\varphi_\gamma\left(V\right)= \sigma_\gamma^{-2}\left( g\left(V\right) - g\left(V_\gamma\right)\right)
= \left(\sigma_\gamma^{-2} - 1\right) g\left(V\right)
= \left(\sigma_\gamma^{-2} - 1\right) \sum_{k=1}^K \sum_{\underset{\ell\neq k}{\ell=1}}^K \tau_k^*V_{k\ell}\tau_\ell .
\]
Thus
\begin{eqnarray*}
Z_\gamma &=& \sigma_\gamma^{-2} \kappa_3\textrm{\large \textbf{1}}_{E_\gamma} (X)\, \sum_{k=1}^K \sum_{\underset{\ell\neq k}{\ell=1}}^K \bigg\{- \frac{1}{2}\bigg(\tau_k^* \left(X_k\otimes X_k\right) V_{k\ell}\tau_\ell +  \tau_\ell^*V_{\ell k} \left(X_k\otimes X_k\right)\tau_k\bigg)\\
& & +	\tau_k^*\left(X_\ell\otimes X_k\right) \tau_\ell\bigg\}
+ \left(\sigma_\gamma^{-2} - 1\right) w\left(\Vert V^{-1/2}X\Vert _\mathcal{X} \right) \sum_{k=1}^K \sum_{\underset{\ell\neq k}{\ell=1}}^K\tau_k^*V_{k\ell}\tau_\ell. 
\end{eqnarray*}
\subsection{Proof of Theorem \ref{loicoeff}}
Arguing as the proof of Theroem 3.3 in \cite{nkiet} we see that  $\sqrt{n}\left(\widehat{\Lambda}_n -\Lambda\right)$ converges in distribution, as $n\rightarrow +\infty$, to the random variable   $W_\gamma=P^\ast U_\gamma P$, where  $P=\sum_{\ell=1}^pe_\ell\otimes \beta^{(\ell )}$ and $\{e_\ell\}_{1\leq \ell\leq q}$ is an orthonormal basis of $\mathcal{X}$. Clearly, $W_\gamma$ has a normal distribution with mean $0$ and covariance operator $\Sigma$ equal to that of $P^\ast Z_\gamma P$ and that will now be explicited. We will first specify $\mathbb{E}(P^\ast Z_\gamma P)$. Considering  the random vector   $Y=V^{-1/2}X$ which has an eliptical distribution, we have  from the proof of Theorem 4.2 in \cite{cator} (see p. 2386) that 
\[
\mathbb{E}\left(\textrm{\large \textbf{1}}_{[0,r(\gamma)]} (\Vert Y\Vert_\mathcal{X})Y\otimes Y\right)=\frac{1}{q}\,\mathbb{E}\left(\textrm{\large \textbf{1}}_{[0,r(\gamma)]} (\Vert Y\Vert_\mathcal{X})\Vert Y\Vert_\mathcal{X}^2\right)\,\mathbb{I}.
\] 
Hence
\begin{eqnarray*}
\mathbb{E}\left(\textrm{\large \textbf{1}}_{E_\gamma} (X)\,\,X_\ell\otimes X_k\right)&=&\tau_kV^{1/2}\mathbb{E}\left(\textrm{\large \textbf{1}}_{[0,r(\gamma)]} (\Vert Y\Vert_\mathcal{X})\,Y\otimes Y\right)V^{1/2}\tau_\ell^\ast\\
&=&\frac{1}{q}\,\mathbb{E}\left(\textrm{\large \textbf{1}}_{[0,r(\gamma)]} (\Vert Y\Vert_\mathcal{X})\Vert Y\Vert_\mathcal{X}^2\right)\,\tau_kV\tau_\ell^\ast\\
&=&\frac{1}{q}\,\mathbb{E}\left(\textrm{\large \textbf{1}}_{E_\gamma} (X)\, \Vert V^{-1/2}X\Vert_\mathcal{X}^2\right)\,V_{k\ell}\\
&=&\frac{\mu}{q}\,V_{k\ell}
\end{eqnarray*}
and, therefore,
\begin{eqnarray*}
\mathbb{E}(Z_\gamma) &=& \sigma_\gamma^{-2} \kappa_3 \sum_{k=1}^K \sum_{\underset{\ell\neq k}{\ell=1}}^K \bigg\{- \frac{1}{2}\bigg(\tau_k^* \mathbb{E}\left(\textrm{\large \textbf{1}}_{E_\gamma} (X)\,\,X_\ell\otimes X_k\right) V_{k\ell}\tau_\ell +  \tau_\ell^*V_{\ell k} \mathbb{E}\left(\textrm{\large \textbf{1}}_{E_\gamma} (X)\,\,X_\ell\otimes X_k\right)\tau_k\bigg)\\
& & +	\tau_k^*\mathbb{E}\left(\textrm{\large \textbf{1}}_{E_\gamma} (X)\,\,X_\ell\otimes X_k\right)\tau_\ell\bigg\}
+ \left(\sigma_\gamma^{-2} - 1\right) \mathbb{E}\left(w\left(\Vert V^{-1/2}X\Vert _\mathcal{X} \right) \right)\sum_{k=1}^K \sum_{\underset{\ell\neq k}{\ell=1}}^K\tau_k^*V_{k\ell}\tau_\ell\\
&=&\frac{\sigma_\gamma^{-2}\mu \kappa_3}{q} \sum_{k=1}^K \sum_{\underset{\ell\neq k}{\ell=1}}^K \bigg\{- \frac{1}{2}\bigg(\tau_k^*  V_{k\ell}\tau_\ell +  \tau_\ell^*V_{\ell k} \tau_k\bigg)+	\tau_k^*V_{k\ell}\tau_\ell\bigg\}\\
& &+ \left(\sigma_\gamma^{-2} - 1\right) \mathbb{E}\left(w\left(\Vert V^{-1/2}X\Vert _\mathcal{X} \right) \right)\sum_{k=1}^K \sum_{\underset{\ell\neq k}{\ell=1}}^K\tau_k^*V_{k\ell}\tau_\ell\\
&=& \left(\sigma_\gamma^{-2} - 1\right) \mathbb{E}\left(w\left(\Vert V^{-1/2}X\Vert _\mathcal{X} \right) \right)\sum_{k=1}^K \sum_{\underset{\ell\neq k}{\ell=1}}^K\tau_k^*V_{k\ell}\tau_\ell.
\end{eqnarray*}
On the other hand, we have 
\begin{eqnarray*}
\mathbb{E}\left(w\left(\Vert V^{-1/2}X\Vert _\mathcal{X} \right) \right)
&=&\kappa_1\mathbb{E}\left(\textrm{\large \textbf{1}}_{E_\gamma} (X) \right) +\kappa_2\mathbb{E}\left(\textrm{\large \textbf{1}}_{E_\gamma} (X)\, \Vert V^{-1/2}X\Vert_\mathcal{X}^2\right)+\kappa_4=
\kappa_1\gamma +\kappa_2\mu+\kappa_4
\end{eqnarray*}
and, consequently,
\begin{eqnarray*}
\mathbb{E}(P^\ast Z_\gamma P) &=&P^\ast \mathbb{E}( Z_\gamma )P =\left(\sigma_\gamma^{-2} - 1\right) \left(\kappa_1\gamma +\kappa_2\mu+\kappa_4\right)\sum_{k=1}^K \sum_{\underset{\ell\neq k}{\ell=1}}^KP^\ast \tau_k^*V_{k\ell}\tau_\ell P\\
&=& \left(\sigma_\gamma^{-2} - 1\right) \left(\kappa_1\gamma +\kappa_2\mu+\kappa_4\right)\\
& &\hspace{1cm}\times\sum_{m=1}^q\sum_{r=1}^q\sum_{k=1}^K \sum_{\underset{\ell\neq k}{\ell=1}}^K(\beta^{(r)}\otimes e_r)\tau_k^*V_{k\ell}\tau_\ell (e_m\otimes\beta^{(m)})\\
&=& \left(\sigma_\gamma^{-2} - 1\right) \left(\kappa_1\gamma +\kappa_2\mu+\kappa_4\right)\\
& &\hspace{1cm}\times\sum_{m=1}^q\sum_{r=1}^q\sum_{k=1}^K \sum_{\underset{\ell\neq k}{\ell=1}}^K(\beta^{(r)}\otimes e_r)(e_m\otimes(\tau_k^*V_{k\ell}\tau_\ell \beta^{(m)}))\\
&=& \left(\sigma_\gamma^{-2} - 1\right) \left(\kappa_1\gamma +\kappa_2\mu+\kappa_4\right)\\
& &\hspace{1cm}\times\sum_{m=1}^q\sum_{r=1}^q\sum_{k=1}^K \sum_{\underset{\ell\neq k}{\ell=1}}^K<\tau_k\beta^{(r)}, V_{k\ell}\tau_\ell \beta^{(m)}>_k\,(e_m\otimes e_r).
\end{eqnarray*}
Note that
\[
Z_\gamma=\sigma_\gamma^{-2} \kappa_3\textrm{\large \textbf{1}}_{E_\gamma} (X)\,Z^{(0)}
+ \left(\sigma_\gamma^{-2} - 1\right) w\left(\Vert V^{-1/2}X\Vert _\mathcal{X} \right) \sum_{k=1}^K \sum_{\underset{\ell\neq k}{\ell=1}}^K\tau_k^*V_{k\ell}\tau_\ell
\]
where
\[
Z^{(0)}=\sum_{k=1}^K \sum_{\underset{\ell\neq k}{\ell=1}}^K \bigg\{- \frac{1}{2}\bigg(\tau_k^* \left(X_k\otimes X_k\right) V_{k\ell}\tau_\ell +  \tau_\ell^*V_{\ell k} \left(X_k\otimes X_k\right)\tau_k\bigg)
+	\tau_k^*\left(X_\ell\otimes X_k\right) \tau_\ell\bigg\};
\]
it is known (see \cite{nkiet}) that
\begin{eqnarray*}
P^\ast Z^{(0)}P=\sum_{m=1}^p\sum_{r=1}^p
\bigg[\sum_{k=1}^K\sum_{\stackrel{\ell=1}{\ell\neq k}}^K &-& \frac{1}{2}\bigg(<\tau_\ell\beta^{(m)},V_{\ell k}X_k>_\ell<\tau_k\beta^{(r)},X_k>_k\\
&+&<\tau_\ell\beta^{(r)},V_{\ell k}X_k>_\ell <\tau_k\beta^{(m)},X_k>_k\bigg)\\
& +&<\tau_\ell\beta^{(m)},X_\ell>_\ell<\tau_k\beta^{(r)},X_k>_k\bigg]\, e_m\otimes e_r.
\end{eqnarray*}
Thus,
\begin{eqnarray*}
P^\ast Z_\gamma P- \mathbb{E}(P^\ast Z_\gamma P) &=&\sigma_\gamma^{-2} \kappa_3\textrm{\large \textbf{1}}_{E_\gamma} (X)\,P^\ast Z^{(0)}P\\
& &- \left(\sigma_\gamma^{-2} - 1\right)\left( w\left(\Vert V^{-1/2}X\Vert _\mathcal{X} \right)-\kappa_1\gamma -\kappa_2\mu-\kappa_4 \right)\sum_{k=1}^K \sum_{\underset{\ell\neq k}{\ell=1}}^KP^\ast\tau_k^*V_{k\ell}\tau_\ell P\\
&=&\sum_{m=1}^p\sum_{r=1}^p
\sum_{k=1}^K\sum_{\stackrel{\ell=1}{\ell\neq k}}^K\sigma_\gamma^{-2} \kappa_3\textrm{\large \textbf{1}}_{E_\gamma} (X)\,\\
& &\times\bigg[ - \frac{1}{2}\bigg(<\tau_\ell\beta^{(m)},V_{\ell k}X_k>_\ell<\tau_k\beta^{(r)},X_k>_k\\
& &\hspace{1.5cm}+<\tau_\ell\beta^{(r)},V_{\ell k}X_k>_\ell <\tau_k\beta^{(m)},X_k>_k\bigg)\\
& &\hspace{2cm}+<\tau_\ell\beta^{(m)},X_\ell>_\ell<\tau_k\beta^{(r)},X_k>_k\bigg]\\
& & -\left(\sigma_\gamma^{-2} - 1\right)\left( w\left(\Vert V^{-1/2}X\Vert _\mathcal{X} \right)-\kappa_1\gamma -\kappa_2\mu-\kappa_4 \right)  \\
& &\hspace{1cm}\times<\tau_k\beta^{(r)}, V_{k\ell}\tau_\ell \beta^{(m)}>_k\,(e_m\otimes e_r)\\
&=&\sum_{m=1}^p\sum_{r=1}^p\mathcal{Y}_{m,r}\,\,(e_m\otimes e_r).
\end{eqnarray*}
Finally
\begin{eqnarray*}
\Sigma&=&\mathbb{E}\bigg(\left(P^\ast Z_\gamma P- \mathbb{E}(P^\ast Z_\gamma P)\right)\widetilde{\otimes}\left(P^\ast Z_\gamma P- \mathbb{E}(P^\ast Z_\gamma P)\right)\bigg)\\
&=&\sum_{1\leq m,r,u,t\leq q}\mathbb{E}\left(\mathcal{Y}_{m,r}\,\mathcal{Y}_{u,t}\right)\,\,(e_m\otimes e_r)\widetilde{\otimes}(e_u\otimes e_t).
\end{eqnarray*}
\subsection{Proof of Theorem \ref{test}}
Under $\mathcal H_0$  we have  $T= 0$ and, therefore,  $\sqrt{n}\widetilde{T}_{n} = \sqrt{n}\left(\widetilde{T}_{n} - T\right)$.  Consequently, from Theorem \ref{loilim} we deduce that  $\sqrt{n}\widetilde{T}_{n}$ converges in distribution,  as  $n \rightarrow +\infty$,  to   $U_\gamma$  which has  a normal distribution in  $\mathcal{L}(\mathcal{X})$   with mean 0 and covariance operator equal to that of  $Z_\gamma$. Since the map \ $A  \mapsto\sum_{k=2}^{K}\sum_{\ell=1}^{k-1}\textrm{tr}\left(\pi_{k\ell}\left(A\right)\pi_{k\ell}\left(A\right)^*\right)  $ is continuous,  we deduce that $ n\widetilde{S}_{n}$ converges in distribution, as  $n \rightarrow +\infty$,
to  
\[
\mathcal{Q}_\gamma= \sum_{k=2}^{K}\sum_{\ell=1}^{k-1}\textrm{tr}\left(\pi_{k\ell}\left( U_\gamma\right)\pi_{k\ell}\left( U_\gamma \right)^*\right).
\]
By  a similar reasoning to that of the proof of  Theorem 4.1 in \cite{nkiet} we obtain that   $\mathcal Q_\gamma = \mathbb W_\gamma^T\mathbb W_\gamma$ where  $\mathbb W_\gamma$ is a random variable having  normal distribution $N(0,\Theta)$  in $\mathbb R^d$  with 
\begin{equation*}\label{cov}
\Theta=\left(
\begin{array}
[c]{ccccc}%
\Theta_{21,21} & \Theta_{21,31} & \Theta_{21,32}  &\cdots & \Theta_{21,KK-1}\\
\Theta_{31,21} & \Theta_{31,31} & \Theta_{31,32}&\cdots & \Theta_{31,KK-1}\\
\vdots & \vdots & \vdots & \cdots&\vdots\\
\Theta_{KK-1,21} & \Theta_{KK-1,31} & \Theta_{KK-1,32}&\cdots & \Theta_{KK-1,KK-1}
\end{array}
\right) ,
\end{equation*}
where $\Theta_{k\ell,rs}$ is the $p_kp_\ell\times p_rp_s$ matrix given by
\begin{equation*}
\Theta_{k\ell,rs}=\left(
\begin{array}
[c]{ccccccccc}%
\gamma_{1111}^{k\ell,rs} & \gamma_{1121}^{k\ell,rs} & \cdots & \gamma_{11p_r1}^{k\ell,rs}&\cdots &\gamma_{111p_s}^{k\ell,rs} & \gamma_{112p_s}^{k\ell,rs} & \cdots & \gamma_{11p_rp_s}^{k\ell,rs}\\
\gamma_{2111}^{k\ell,rs} & \gamma_{2121}^{k\ell,rs} & \cdots & \gamma_{21p_r1}^{k\ell,rs}&\cdots &\gamma_{211p_s}^{k\ell,rs} & \gamma_{212p_s}^{k\ell,rs} & \cdots & \gamma_{21p_rp_s}^{k\ell,rs}\\
\vdots &\vdots &\cdots &\vdots &\cdots &\vdots & \vdots &\cdots &\vdots \\
\gamma_{p_k111}^{k\ell,rs} & \gamma_{p_k121}^{k\ell,rs} & \cdots & \gamma_{p_k1p_r1}^{k\ell,rs}&\cdots &\gamma_{p_k11p_s}^{k\ell,rs} & \gamma_{p_k12p_s}^{k\ell,rs} & \cdots & \gamma_{p_k1p_rp_s}^{k\ell,rs}\\
\vdots &\vdots &\cdots &\vdots &\cdots &\vdots & \vdots &\cdots &\vdots \\
\gamma_{1p_\ell 11}^{k\ell ,rs} & \gamma_{1p_\ell 21}^{k\ell ,rs} & \cdots & \gamma_{1p_\ell p_r1}^{k\ell ,rs}&\cdots &\gamma_{1p_\ell 1p_s}^{k\ell ,rs} & \gamma_{1p_\ell 2p_s}^{k\ell ,rs} & \cdots & \gamma_{1p_\ell p_rp_s}^{k\ell ,rs}\\
\gamma_{2p_\ell 11}^{k\ell,rs} & \gamma_{2p_\ell 21}^{k\ell ,rs} & \cdots & \gamma_{2p_\ell p_r1}^{k\ell,rs}&\cdots &\gamma_{2p_\ell 1p_s}^{k\ell,rs} & \gamma_{2p_\ell 2p_s}^{k\ell,rs} & \cdots & \gamma_{2p_\ell p_rp_s}^{k\ell,rs}\\
\vdots &\vdots &\cdots &\vdots &\cdots &\vdots & \vdots &\cdots &\vdots \\
\gamma_{p_kp_\ell 11}^{k\ell,rs} & \gamma_{p_kp_\ell 21}^{k\ell,rs} & \cdots & \gamma_{p_kp_\ell p_r1}^{k\ell,rs}&\cdots &\gamma_{p_kp_\ell 1p_s}^{k\ell,rs} & \gamma_{p_kp_\ell 2p_s}^{k\ell,rs} & \cdots & \gamma_{p_kp_\ell p_rp_s}^{k\ell,rs}
\end{array}
\right) ,
\end{equation*}
and        
\begin{eqnarray*}
\gamma_{ijpt}^{kl,ru} &=& < \mathbb{E}\bigg(\pi_{k\ell}\left( U_\gamma\right)\widetilde{\otimes}\pi_{rs}\left(U_\gamma \right)\bigg)\left(e_j^{(l)}\otimes e_i^{(k)}\right) , e_t^{(u)}\otimes e_p^{(r)}>\\
&=&< \mathbb{E}\bigg(\pi_{k\ell}\left( Z_\gamma\right)\widetilde{\otimes}\pi_{rs}\left(Z_\gamma \right)\bigg)\left(e_j^{(l)}\otimes e_i^{(k)}\right) , e_t^{(u)}\otimes e_p^{(r)}>
\end{eqnarray*}
where   $\widetilde{\otimes}$  denotes the tensor product related to the inner product of operators
$<  A , B > = tr (AB^*)$ and $\{e_i^{(k)}\}_{1\leq i\leq p_k}$ is an orthonormal basis of $\mathcal{X}_k$. Since under $\mathcal H_0$,   $Z_\gamma$ becomes   
\begin{eqnarray*}
Z_{\gamma}&=&\sigma_\gamma^{-2}\kappa_3\textrm{\large \textbf{1}}_{E_\gamma} (X) \sum_{k=1}^K \sum_{\underset{\ell\neq k}{\ell=1}}^K \tau_k^*\left(X_\ell\otimes X_k\right) \tau_\ell
\end{eqnarray*}
we obtain
\begin{eqnarray*}
\pi_{k\ell}\left( Z_\gamma\right)&=&\sigma_\gamma^{-2}\kappa_3\textrm{\large \textbf{1}}_{E_\gamma} (X) \sum_{j=1}^K \sum_{\underset{m\neq j}{m=1}}^K \tau_k\tau_j^*\left(X_m\otimes X_j\right) \tau_m\tau_{\ell}^\ast\\
&=&\sigma_\gamma^{-2}\kappa_3\textrm{\large \textbf{1}}_{E_\gamma} (X) \sum_{j=1}^K \sum_{\underset{m\neq j}{m=1}}^K\delta_{kj}\,\delta_{\ell m} \left(X_m\otimes X_j\right) \\
&=&\sigma_\gamma^{-2}\kappa_3\textrm{\large \textbf{1}}_{E_\gamma} (X) \left(X_\ell\otimes X_k\right).
\end{eqnarray*}
Hence
\begin{eqnarray*}
\gamma_{ijpt}^{kl,ru} &=&\sigma_\gamma^{-4} \kappa_3^2 \mathbb{E}\bigg(\textrm{\large \textbf{1}}_{E_\gamma} (X)<\left(\left( X_\ell\otimes X_k \right)\widetilde{\otimes}\left( X_s\otimes X_r  \right)\right)\left(e_j^{(\ell)}\otimes e_i^{(k)}\right) , e_t^{(u)}\otimes e_p^{(r)}>\bigg)\\
&=&\sigma_\gamma^{-4} \kappa_3^2\mathbb{E}\bigg(\textrm{\large \textbf{1}}_{E_\gamma} (X)<  X_l\otimes X_k , e_j^{(l)}\otimes e_i^{(k)}>< X_s\otimes X_r , e_t^{(u)}\otimes e_p^{(r)}>\bigg)\\
&=&\sigma_\gamma^{-4} \kappa_3^2\mathbb{E}\bigg(\textrm{\large \textbf{1}}_{E_\gamma} (X)< X_k , e_i^{(k)}>_k<  X_r , e_p^{(r)}>_r <  X_\ell , e_j^{(l)}>_\ell< X_s , e_t^{(u)}>_s\bigg).
\end{eqnarray*}
Note that  under $\mathcal{H}_0$ we have $V=\mathbb{I}$, then $X$ has a spherical distribution with density $f_X(x)=h(\Vert x\Vert^2_\mathcal{X})$. Therefore, if $(k,\ell)=(r,u)$ and $(i,j)=(p,t)$ with $\ell\neq k$ and $u\neq r$, then $\gamma_{ijpq}^{kl,ru} $ equals an integral  of the form
\[
\int z\left (\Vert x\Vert^2_\mathcal{X}\right)\,x^2_a\,x^2_b\,\,dx,
\]
with $a\neq b$, where $z$ is a suitable function from $[0,+\infty[$ to itself.  Then from Lemma 1 in \cite{lopuhaa}, we deduce that
\[
\int z\left (\Vert x\Vert^2_\mathcal{X}\right)\,x^2_a\,x^2_b\,\,dx=\frac{1}{q(q+1)}\int z\left (\Vert x\Vert^2_\mathcal{X}\right)\,\Vert x\Vert^4_\mathcal{X}\,dx
\]
and, therefore, that 
\[
\gamma_{ijpt}^{kl,ru}=\frac{\sigma_\gamma^{-4} \kappa_3^2}{q(q+1)}\mathbb{E}\bigg(\textrm{\large \textbf{1}}_{E_\gamma} (X)\Vert  X\Vert^4_\mathcal{X}\bigg).
\]
Otherwise, if  one of the conditions $(k,\ell)=(r,u)$, $(i,j)=(p,t)$,  $\ell\neq k$,  $u\neq r$ does not hold then  $\gamma_{ijpt}^{kl,ru} $ equals an integral  of the form
\[
\int z\left (\Vert x\Vert^2_\mathcal{X}\right)\,x_a\,x_b\,x_c\,x_d\,\,dx
\]
in which at least two of the indices $a$, $b$, $c$, $d$ are different. From elementary calculus obtained by changing to spherical coordinates we obtain that this integral equals $0$ and, therefore, $\gamma_{ijpt}^{kl,ru} =0$. We deduce that 
\begin{eqnarray*}
\Theta &=& \frac{\sigma_\gamma^{-4} \kappa_3^2}{q(q+1)}\mathbb{E}\bigg(\textrm{\large \textbf{1}}_{E_\gamma} (X)\Vert  X\Vert^4_\mathcal{X}\bigg)\, I_d
\end{eqnarray*}
where $I_d$ is the $d\times d$ identity matrix. Thus,  $\mathcal{Q}_\gamma = \frac{\sigma_\gamma^{-4} \kappa_3^2}{q(q+1)}\mathbb{E}\bigg(\textrm{\large \textbf{1}}_{E_\gamma} (X)\Vert  X\Vert^4_\mathcal{X}\bigg)\mathcal{Q}^\prime$  where $\mathcal{Q}^\prime$ is a random variable with distribution equal to  $\chi_d^2$.
\subsection{Proof of Theorem \ref{if2}}
We have:
\begin{eqnarray*}
\frac{\partial^2 \mathbb{S}_\gamma(\mathbb{P}_{\varepsilon,x})}{\partial^2\varepsilon}
&=&  \sum_{k=2}^{K} \sum_{\ell=1}^{k - 1}\bigg\{\textrm{tr}\bigg[\pi_{k\ell}\left(\frac{\partial^2 \mathbb{T}_{1,\gamma}(\mathbb{P}_{\varepsilon,x})}{\partial\varepsilon^2}\right)\bigg(\pi_{k\ell}\left( \mathbb{T}_{1,\gamma}(\mathbb{P}_{\varepsilon,x})\right)\bigg)^*\bigg]\\
& & \hspace{1.5cm}+ 2\textrm{tr}\bigg[\pi_{k\ell}\left(\frac{\partial \mathbb{T}_{1,\gamma}(\mathbb{P}_{\varepsilon,x})}{\partial\varepsilon}\right)\bigg(\pi_{k\ell}\left(\frac{\partial \mathbb{T}_{1,\gamma}(\mathbb{P}_{\varepsilon,x})}{\partial\varepsilon}\right)\bigg)^*\bigg]\\
& &\hspace{1.5cm}+ \textrm{tr}\bigg[ \pi_{k\ell}\left( \mathbb{T}_{1,\gamma}(\mathbb{P}_{\varepsilon,x})\right)\bigg(\pi_{k\ell}\left(\frac{\partial^2 \mathbb{T}_{1,\gamma}(\mathbb{P}_{\varepsilon,x})}{\partial\varepsilon^2}\right)\bigg)^*\bigg]\bigg\}
\end{eqnarray*}
and, therefore, 
\begin{eqnarray*}
\textrm{IF}^{(2)}(x; S_\gamma, \mathbb{P}_X)&=& \frac{\partial^2 \mathbb S_\gamma(\mathbb{P}_{,X})}{\partial\varepsilon^2}\bigg\vert_{\varepsilon =0}\\
&=&  \sum_{k=2}^{K} \sum_{\ell=1}^{k - 1}\bigg\{\textrm{tr}\bigg[ \pi_{k\ell}\left(\frac{\partial^2 \mathbb{T}_{1,\gamma}(\mathbb{P}_{\varepsilon,x})}{\partial\varepsilon^2}\bigg\vert_{\varepsilon =0}\right)\pi_{k\ell}\left( T\right)^*\bigg]\\
& & \hspace{1.5cm}+ 2\textrm{tr}\bigg[\pi_{k\ell}\left(\textrm{IF}(x; T_\gamma, \mathbb{P}_X)\right)\pi_{k\ell}\left(\textrm{IF}(x; T_\gamma, \mathbb{P}_X)\right)^*\bigg]\\
& &\hspace{1.5cm}+ \textrm{tr}\bigg[ \pi_{k\ell}\left( T\right)\bigg(\pi_{k\ell}\left(\frac{\partial^2 \mathbb{T}_{1,\gamma}(\mathbb{P}_{\varepsilon,x})}{\partial\varepsilon^2}\bigg\vert_{\varepsilon =0}\right)\bigg)^*\bigg]\bigg\}.
\end{eqnarray*}
Since, under $\mathcal{H}_0$, one has $T=0$ and 
\begin{eqnarray*}
\pi_{k\ell}\left(\textrm{IF}( x ; T_\gamma, \mathbb{P}_X)\right)&=&-\frac{\sigma_\gamma^{-2}}{2\kappa_0}\textrm{\large \textbf{1}}_{E_\gamma} (x)\,\sum_{j=1}^K \sum_{\underset{m\neq j}{m=1}}^K \pi_{k\ell}\left( \tau_j^*\left(x_m\otimes x_j\right)\tau_m\right)\\
&= &-\frac{\sigma_\gamma^{-2}}{2\kappa_0}\textrm{\large \textbf{1}}_{E_\gamma} (x)\,\sum_{j=1}^K \sum_{\underset{m\neq j}{m=1}}^K \tau_k  \tau_j^*\left(x_m\otimes x_j\right)\tau_m\tau_\ell^\ast\\
&= &-\frac{\sigma_\gamma^{-2}}{2\kappa_0}\textrm{\large \textbf{1}}_{E_\gamma} (x)\,\sum_{j=1}^K \sum_{\underset{m\neq j}{m=1}}^K \delta_{jk}\delta_{m\ell} \left(x_m\otimes x_j\right) \\
&=&-\frac{\sigma_\gamma^{-2}}{2\kappa_0}\textrm{\large \textbf{1}}_{E_\gamma} (x)\, x_\ell\otimes x_k
\end{eqnarray*}
 it follows
\begin{eqnarray*}
\textrm{IF}^{(2)}(x; S_\gamma, \mathbb{P}_X)&=& 2 \sum_{k=2}^{K} \sum_{\ell=1}^{k - 1}\textrm{tr}\bigg[\pi_{k\ell}\left(\textrm{IF}(x; T_\gamma, \mathbb{P}_X)\right)\pi_{k\ell}\left(\textrm{IF}(x; T_\gamma, \mathbb{P}_X)\right)^*\bigg]\\
&=&\frac{\sigma_\gamma^{-4}}{4\kappa_0^2}\textrm{\large \textbf{1}}_{E_\gamma} (x) \,\sum_{k=2}^{K} \sum_{\ell=1}^{k - 1}\textrm{tr}\bigg[\left(x_\ell\otimes x_k\right)\left(x_\ell\otimes x_k\right)^*\bigg].
\end{eqnarray*}
Then, using $(u\otimes v)^\ast=v\otimes u$, $(x\otimes y)(z\otimes t)=<x,t>\,(z\otimes y)$ and $\textrm{tr}(x\otimes y)=<x,y>$ (see \cite{dauxois}), we finally obtain
\[
\textrm{IF}^{(2)}(x; S_\gamma, \mathbb{P}_X)=\frac{\sigma_\gamma^{-4}}{4\kappa_0^2}\textrm{\large \textbf{1}}_{E_\gamma} (x)\sum_{k=2}^{K} \sum_{\ell=1}^{k - 1}\Vert x_k\Vert  _k^2\,\,\Vert x_\ell\Vert  _\ell^2.
\]

\end{document}